\documentclass[12pt,a4paper]{article}
\usepackage{mathrsfs}
\usepackage{latexsym}
\usepackage{amssymb}
\usepackage{graphicx}
\usepackage{amsmath}
\usepackage[usenames]{color}
\newtheorem{theorem}{Theorem}
\newtheorem{define}[theorem]{Definition}

\newtheorem{coro}[theorem]{Corollary}
\newtheorem{lemma}[theorem]{Lemma}
\newtheorem{remark}[theorem]{Remark}
\newtheorem{example}[theorem]{Example}
\newenvironment{pf}{\vspace{3mm}\indent{\bf  Proof. }}{\hfill $\Box$ \vspace{3mm}}

\begin{document}
\title{Moran-Type Iterated Function Systems and Dimensions of Moran Self-Similar Sets\author{Yong-Shen Cao, Qi-Rong Deng and Ming-Tian Li}
\thanks{The work is supported by NNSFC (No. 11971109), the Program for Probability and Statistics: Theory and Application (No. IRTL1704),  the program for innovative research team in science and technology in Fujian province university (No. IRTSTFJ) and  the Natural Science Foundation
of Fujian Province (No. 2023J01298).}}

\maketitle

\begin{abstract} 
Moran-type iterated function systems (Moran-type IFS or MIFS) are defined by a sequence of  iterated function systems, and their basic theoretical framework is established. We define Moran-type attractors and invariant probability measures associated with a sequence of probability weight vectors. Furthermore, separation conditions for MIFS are introduced, and the dimension theory of Moran-type self-similar sets is investigated. Appropriate examples are provided to illustrate and support the definitions and results.

\medskip

\end{abstract}

\section{Introduction}
\setcounter{equation}{0}\setcounter{theorem}{0}
Since B. B. Mandelbrot proposed the concept of fractal in the 1980s (see \cite{[1]}), fractal geometry, as a new branch of mathematics, has  progressed significantly during the past 45 years.

\medskip

Moran set, as a generalization of self similar sets, is inspired by the idea of P. A. Moran (see \cite{[4]}). The constructions are direct generalizations of the construction of the middle third Cantor set. Z. Y. Wen et al. called such sets as Moran sets (see, for example, \cite{[5]}, \cite{[6]}, \cite{[7]}, \cite{[8]}, \cite{[9]}, \cite{[10]}, \cite{[11]}, \cite{[12]}, \cite{[DWXX]} and the references therein). In these literatures, the authors mainly considered a special class of Moran sets by the following construction: First, fix a structurally simple compact set $E_1$ (an interval, a square, etc.). Then take several non overlapping subsets of $E_1$ that are similar with $E_1$ to generate a subset $E_2\subset E_1$. Then take several non-overlapping subsets of $E_2$ that are similar with $E_1$ to generate a subset $E_3\subset E_2$. Continuing this process generates a limit set, i.e., the so-called Moran set. The dimension problem of Moran sets in different situations are studied. These results indicate that when the contraction ratios and translations of these similar maps have different patterns, the corresponding dimensions and Hausdorff measures of Moran sets have different patterns and may not be $s$-sets. Some further results on Moran sets are also obtained in \cite{[13]}, \cite{[DMW]}, \cite{[LLMX]}, etc.

\medskip

Since J. Hutchinson introduced the powerful tool of iterated function systems in \cite{[2]}, it plays an extremely important role in the field of fractal geometry as an important tool for defining fractal sets and fractal measures and, has attracted great interest of a large number of mathematicians. Under the framework established by Hutchinson, a large number of research achievements have been made. The number of citations of reference \cite{[2]} indicates its significant impact. Just as the well known fractal geometer M. Barnsley (see \cite{[3]}) stated: iterated function systems have been at the heart of fractal geometry almost from its origins. Based on this observation, we can see that applying Hutchinson's ideas to the construction of Moran sets should be a very meaningful work. This is our first motivation in this article.

\medskip

Our second motivation in this article comes from the investigation on the spectrality of fractal measures. Since P. Jorgenson and S. Pederson found the first example of singular spectral probability measure, the spectrality  and non-spectrality of self-affine measures have attracted great interest of a large number of mathematicians, a large number of research achievements have been made during the past twenty years (see, for example,  \cite{[JP]}, \cite{[JP1]}, \cite{[S]}, \cite{[S2]}, \cite{[HL]}, \cite{[D1]}, \cite{[Dai2]}, \cite{[DHL]}, \cite{[DHLu]}, \cite{AH}, \cite{[28]}, \cite{[27]}, \cite{[DeL]}, \cite{[DDP]}, \cite{[DJ1]}, \cite{[DJ2]}, \cite{[22]}, \cite{[32]}, \cite{[FHL]}, \cite{[21]}, \cite{[HL]}, \cite{[LW2]}, \cite{[30]}, \cite{[31]}, \cite{[Li1]}, \cite{[Li2]}, \cite{[Li3]}, \cite{[Li4]}, \cite{[26]}, \cite{[LDL]}, \cite{[29]}, \cite{[25]}, \cite{[23]} and the references therein). As a generalization of spectrality of self-affine measures, L. X. An and X. G. He first studied the spectrality of a class of Moran-type self-affine measures (can be written as infinite convolutions of uniform distributions on finite sets), then this new field immediately attracted the interest of a large number of mathematicians, many important results have been made during the past 10 years (see, for example, \cite{[17]}, \cite{[18]}, \cite{[19]}, \cite{[20]}, \cite{[21]}, \cite{[22]}, \cite{[23]}, \cite{[24]}, \cite{[26]}, \cite{[29]}, \cite{[33]} and the references therein).

\begin{define}\label{def1.1}
Let $X\subset{\mathbb R}^d$ be a compact set, $\Phi_n=\left\{\phi_{n,j}\right\}_{j=1}^{N_n}$ with $2\le N_n<+\infty$ be a collection of bi-Lipschitz maps from $X$ to itself. If there exist constants $0<c_{1,n}\le c_{2,n}<1$ and a metric $\rho(\cdot,\ \cdot)$ on ${\mathbb R}^d$ such that
\begin{equation}\label{eq1.1}\left\{\begin{array}{l}
\phi_{n,j}(X)\subseteq X,\ \ \ \lim\limits_{n\to+\infty}\prod_{j=1}^nc_{2,j}=0\\ \ c_{1,n}\rho(x,\ y)\le\rho(\phi_{n,j}(x),\ \phi_{n,j}(y))\le {}{c_{2,n}}\rho(x,\ y),\ \forall\ x,\ y\in X
\end{array}\right.
\end{equation}
whenever $1\le j\le N_n$ and $n\ge 1$. We call the sequence of bi-Lipschitz iterated function systems $\left\{\Phi_n\right\}_{n=1}^\infty$ a Moran-type iterated function system (MIFS in short) on $(X,\ \rho(\cdot,\ \cdot))$.
\end{define}

\begin{remark}\label{remark1.2} By Definition \ref{def1.1}, we see that $\left\{\Phi_n\right\}_{n=k}^\infty$ is also a Moran-type iterated function system for any $k>1$.
\end{remark}

\begin{remark}\label{remark1.3} Let $\left\{\phi_{j}\right\}_{j=1}^{m}$ be an IFS in the compact metric space $(X,\ \rho(\cdot,\ \cdot))$ as in \cite{[2]}. By defining $\phi_{n,j}=\phi_{j}$ and $N_n=m$ for all $n\ge1$, we see that $\left\{\Phi_n\right\}_{n=1}^\infty$ is also a Moran-type iterated function system. {{Hence, the concept of MIFS is a generalization of  IFS in \cite{[2]}}}.
\end{remark}

\medskip

The paper is organized as follows. In Section 2, the definitions of attractors and invariant measures are given after some fundamental theorems having been proven. In Section 3, the open set condition, the weak separation condition, the bounded distortion property for iterated function systems are generalized to the case of Moran-type iterated function systems. The relationship between these separation conditions for iterated function systems and their generalizations for Moran-type iterated function systems are also discussed. In Section 4, dimensions of Moran-type self-similar sets are investigated. Finally, some examples are given in Section 5 to explain why we give these definitions. Furthermore, some counter examples are also given to show that our assumptions in theorems are needed.

\medskip

The following notations will be used in the following sections.

\medskip

{\bf\large Notations.} For any subset $E\subseteq{\mathbb R}^d$, we let $\dim_{\rm H}(E)$, $\dim_{\rm P}(E)$, $\dim_{\rm B}(E)$,
$\mathcal{H}^s(E)$, $\mathcal{L}^d(E)$, $|E|$ and $E^\circ$
denote the Hausdorff dimension, the packing dimension, the box dimension, the
$s$-dimensional Hausdorff measure, the $d$-dimensional Lebesgue measure, the Euclidean diameter
and the interior of $E$, respectively.

\medskip

Given an MIFS $\left\{\Phi_n\right\}_{n=1}^\infty$ on compact metric space $(X,\ \rho)$ as in Definition \ref{def1.1}, we call
$\Sigma_n=\left\{1,2,\cdots,N_n\right\}( n\geq 1)$ the basic symbolic space.  We write $\Sigma_n^k=\Sigma_n\times\cdots \times \Sigma_{n+k-1}$ for the product of  $k$ consecutive basic spaces. Every element $J=j_{n}\cdots j_{n+k-1}\in \Sigma_n^k$ is called a $k$-length word starting from the $n$-th layer. Let $\Sigma_n^*$ stand for the set of all finite length words starting from the $n$-th layer, i.e.,
\[ \Sigma_n^*=\bigcup\limits_{k=0}^{+\infty}\Sigma_n^k,\ \ \forall\ n\ge1,
\]
where $\Sigma_n^0=\{\vartheta\}$ and $\vartheta$ is the empty word. Furthermore, we write
\[
\Sigma_n^{\mathbb N}:=\left\{j_nj_{n+1}j_{n+2}\cdots j_{k}\cdots:\ 1\le j_k\le N_k\mbox{ for all }k\ge n\right\}, \quad n\geq1.
\]
Given a finite word $J=j_nj_{n+1}j_{n+2}\cdots j_{n+k-1}\in \Sigma_n^k$, we let  $\Sigma_n^{\mathbb N}(J)$ stand for the cylinder set starting from the $n$-th layer, i.e.,
\[
\Sigma_n^{\mathbb N}(J)=\left\{j_nj_{n+1}j_{n+2}\cdots j_{n+k-1}\sigma:\ \sigma\in\Sigma_{n+k}^{\mathbb N}\right\}.
\]
For any $J=j_{n} j_{n+1}\cdots j_{n+k-1}\in \Sigma_n^k$, we define $\phi_{n,J}=\phi_{n,j_{n}}\circ\phi_{n+1,j_{n+1}}\cdots\circ\phi_{n+k-1,j_{n+k-1}}$ when $k>0$ and $\phi_{n,\vartheta}$ is the identity map on $X$. Accordingly, for the attractors $K_n$ {(see the definition in Theorem $2.1$)} and probability weights $p_{n,j}$, we write $K_{n,J}=\phi_{n,j_{n}}\circ\phi_{n+1,j_{n+1}}\cdots\circ\phi_{n+k-1,j_{n+k-1}}(K_{n+k})$ and $p_{n,J}=p_{n,j_{n}}p_{n+1,j_{n+1}}\cdots p_{n+k-1,j_{n+k-1}}$, respectively.

\medskip

For contraction ratios of maps in the MIFS, we write
\begin{equation}\label{eq1.2}
r_J=\min\left\{\frac{\rho(\phi_{1,J}(x),\ \phi_{1,J}(y))}{\rho(x,\ y)}:\ x\ne y\in X\right\},\quad\ J\in \Sigma_1^*,
\end{equation}
\begin{equation}\label{eq1.3}
R_J=\max\left\{\frac{\rho(\phi_{1,J}(x),\ \phi_{1,J}(y))}{\rho(x,\ y)}:\ x\ne y\in X\right\},\quad\ J\in \Sigma_1^*,
\end{equation}
and
\begin{equation}\label{eq1.4}
\left\{\begin{array}{l}
{\mathcal I}_b=\left\{j_1 j_2 \cdots j_n\in\Sigma_1^*:\ R_{j_1 j_2\cdots j_n} \le b < R_{j_1 j_2\cdots j_{n-1}} \right\},\medskip\\ {\mathcal A}_b=\left\{\phi_{1,J}:\ J\in{\mathcal I}_b\right\},
\end{array}\right.\ \ \forall\ b\in(0,\ 1).
\end{equation}

\medskip

\begin{remark} If the maps $\phi_{n,j}$ are all similarities or conformal maps, we will use the Euclidean norm $\|\cdot\|$ instead of the metric $\rho(\cdot,\ \cdot)$. If the maps $\phi_{n,j}$ are affine maps, then $\phi_{n,j}$ maybe not contractive under the Euclidean norm $\|\cdot\|$.  Thus, an appropriate metric $\rho(\cdot,\ \cdot)$ is desirable since it is possible that $R_{j_1 j_2\cdots j_n} > R_{j_1 j_2\cdots j_{n-1}}$ if we use the Euclidean norm $\|\cdot\|$ in the definition of $R_{j_1 j_2\cdots j_n}$.
\end{remark}

\section{Definitions of attractors and invariant measures}
\setcounter{equation}{0}\setcounter{theorem}{0}

\begin{theorem}\label{th2.1}
Let $\left\{\Phi_n\right\}_{n=1}^\infty$ be an MIFS on $(X,\ \rho)$ as in Definition \ref{def1.1}. We have the following statements.

\medskip

{\rm (i)} The following sets
\begin{equation}\label{eq2.1}
K_n:=\left\{\lim\limits_{k\to+\infty}\phi_{n,j_{n}j_{n+1}\cdots j_{n+k}}(a):\ 1\le j_{n+k}\le N_{n+k},\ k\ge0\right\},\ \forall \ n\ge1
\end{equation}
are well defined compact subsets of $X$ and are independent of $a\in X$. Moreover, the above $\left\{K_n\right\}_{n=1}^\infty$ is the unique sequence of compact subsets of $X$ such that
\begin{equation}\label{eq2.2}
K_n=\bigcup\limits_{j=1}^{N_n} \phi_{n,j}(K_{n+1}),\quad\forall \ n\ge1.
\end{equation}

\medskip

{\rm (ii)} For any given sequence of probability weights ${\mathcal P}_n=\left\{p_{n,1},p_{n,2},\cdots, p_{n,N_n}\right\}$ {\rm (}i.e., $\sum\limits_{j=1}^{N_n}p_{n,j}=1$ and $p_{n,j}>0$ for all $n>0$ and $1\le j\le N_n${\rm )}, there exists a unique sequence of probability measures $\left\{\mu_n\right\}_{n=1}^\infty$ supported on $\left\{K_n\right\}_{n=1}^\infty$, respectively, such that
\begin{equation}\label{eq2.3}
\mu_n =\sum\limits_{j=1}^{N_n} p_{n,j}\mu_{n+1} \circ\phi_{n,j}^{-1},\quad n\ge1.
\end{equation}
\end{theorem}

\begin{pf} (i) We divide the proof of conclusion (i) into three claims.

\medskip

{\bf Claim 1}. Each $K_n$ is well defined and is independent of $a\in X$, $n=1,\ 2,\ \cdots$.

\medskip

Proof of Claim 1. For any $a,\ b\in X,\ n>0$ and  any given sequence $\left\{j_k \right\}_{k=n}^{+\infty}$ with $1\le j_k\le N_k$, since $(X,\ \rho(\cdot,\ \cdot))$ is a compact metric space, $\left\{\phi_{n,j_nj_{n+1}\cdots j_k}(a)\right\}_{k=n}^\infty$ contains at least one convergent subsequence. Hence, there is a sequence of positive integers $\{k_m\}_{m=1}^\infty$ and a point $x_0\in X$ such that
\begin{equation}\label{eq2.4}
\lim_{m\to\infty}\phi_{n,j_nj_{n+1}\cdots j_{k_m}}(a)=x_0.
\end{equation}
For any positive integer $p>k_1$, we have $k_m\le p<k_{m+1}$ for some $m\ge1$. The assumption $\phi_{n,j}(X)\subseteq X$ implies
\[
\rho(\phi_{n,j_nj_{n+1}\cdots j_{k_m}}(a),\ \phi_{n,j_nj_{n+1}\cdots j_{p}}(b))\le |\phi_{n,j_nj_{n+1}\cdots j_{n_m}}(X)|_\rho.
\]
By noting $\phi_{n,j_nj_{n+1}\cdots j_{k_m}}(a),\ \phi_{n,j_nj_{n+1}\cdots j_{p}}(b)\in\phi_{n,j_nj_{n+1}\cdots j_{k_m}}(X)$, we have
\[\begin{array}{rl}
&\rho(x_0,\ \phi_{n,j_nj_{n+1}\cdots j_{p}}(b))\\
\le& \rho(x_0,\ \phi_{n,j_nj_{n+1}\cdots j_{k_m}}(a))+ \rho(\phi_{n,j_nj_{n+1}\cdots j_{k_m}}(a),\ \phi_{n,j_nj_{n+1}\cdots j_{p}}(b))\\
\le& \rho(x_0,\ \phi_{n,j_nj_{n+1}\cdots j_{k_m}}(a))+ |\phi_{n,j_nj_{n+1}\cdots j_{k_m}}(X)|_\rho.
\end{array}
\]
The assumption \eqref{eq1.1} implies $|\phi_{n,j_nj_{n+1}\cdots j_{n_m}}(X)|_\rho\le \prod_{j=n}^{n_m}c_{2,j}|X|_\rho  \to 0$ as $m\to+\infty$. Combining this inequality and \eqref{eq2.4} implies
\[
\lim_{p\to\infty}\rho(x_0,\ \phi_{n,j_nj_{n+1}\cdots j_{p}}(b))=0,
\]
so $\lim\limits_{p\to\infty}\phi_{n,j_nj_{n+1}\cdots j_{p}}(b)=x_0$ for all $b\in X$. Then, the uniqueness of $\lim\limits_{p\to\infty}\phi_{n,j_nj_{n+1}\cdots j_{p}}(b)$ shows that $\lim\limits_{p\to\infty}\phi_{n,j_nj_{n+1}\cdots j_{p}}(a)=x_0$ and the limit is independent of $a\in X$.
Therefore, $K_n$ is well defined and  independent of $a\in X$ for all $n>0$.

\medskip

{\bf Claim 2}. Each $K_n$ is compact.

\medskip

Proof of Claim 2. We need only to prove that $K_1$ is compact. Let $\left\{x_n \right\}_{n\geq1}$ be a Cauchy sequence of $K_1$. Claim 1 implies that there are sequences $\left\{j_{n ,i}\right\}_{i=1}^{+\infty}\in\Sigma_1^{\mathbb N}$ such that
\[
x_n =\lim_{i\to+\infty}\phi_{1,j_{n ,1}}\circ\phi_{2,j_{n,2}}\circ\cdots\circ\phi_{i,j_{n,i}}(a).
\]
Noting two facts that $1\le j_{n ,1}\le N_1$, and  $\bigcup\limits_{\tau=1}^{N_1}\left\{n:\ j_{n,1}=\tau\right\}=\mathbb{N}$, we see that there is a positive integer $\tau_1$ such that $\left\{n:\ j_{n,1}=\tau_1\right\}$ is an infinite set. Similarly, 
there is a positive integer $\tau_2$ such that $\left\{n:\ j_{n,1}=\tau_1,\ j_{n,2}=\tau_2\right\}$ is an infinite set. Continuing in this way, we see that there exists a sequence of positive integers $\{\tau_{p}\}_{p\geq1}$ such that the set $\left\{n :\ j_{n,i}=\tau_i,\ i=1,2,\cdots,p\right\}$ is an infinite set for any $p\ge 1$. Then we can find a sequence of integers $\left\{n_j\right\}_{j\geq1}$ such that $n_1<n_2<\cdots<n_j<\cdots$ and $n_p\in\left\{n:\ j_{n,i}=\tau_i,\ i=1,2,\cdots,p\right\}$ for all $p\ge 1$.  Let
\[
x=\lim_{k\to+\infty}\phi_{1,\tau_1}\circ\phi_{2,\tau_{2}}\circ\cdots\circ\phi_{k,\tau_{k}}(a).
\]
Then, we obtain $x\in K_1$ and
\[
\rho(x,\ x_{n_k})\le |\phi_{1,\tau_1}\circ\phi_{2,\tau_{2}}\circ\cdots\circ\phi_{k,\tau_{k}}(X)|_\rho\to0\mbox{ as }k\to\infty
\]
by using $x,\ x_{n_k}\in \phi_{1,\tau_1}\circ\phi_{2,\tau_{2}}\circ\cdots\circ\phi_{k,\tau_{k}}(X)$ and the assumption \eqref{eq1.1}. Hence, we have $x=\lim\limits_{k\to\infty}x_{n_k}=\lim\limits_{n\to\infty}x_{n}$ by noting that $\left\{x_n\right\}_{n\geq1}$ is a Cauchy sequence of $K_1$. This means that any Cauchy sequence of $K_1$ converges to a point in $K_1$, which implies $K_1$ is compact.

\medskip

{\bf Claim 3}. The $\left\{K_n\right\}_{n\ge1}$ defined in \eqref{eq2.1} is the unique sequence of compact sets satisfying \eqref{eq2.2}.

\medskip

{ Proof of Claim 3}.  It is clear that \eqref{eq2.2} follows from \eqref{eq2.1}, so the $\left\{K_n\right\}_{n\ge1}$ defined in \eqref{eq2.1} is a sequence of compact sets satisfying \eqref{eq2.2}. We then prove the uniqueness in the following.

\medskip

Let $\left\{U_n\right\}_{n\geq1}$ be another sequence of compact sets satisfying \eqref{eq2.2}. For any $x_1\in U_1$, \eqref{eq2.2} shows there is $x_2\in U_2$ and $1\le j_1\le N_1$ such that $x_1=\phi_{1,j_1}(x_2)$. Using \eqref{eq2.2} repeatedly shows that there exist $x_k\in U_k$ and $1\le j_k\le N_k$ such that $x_k=\phi_{k,j_k}(x_{k+1})$ for $k=1,2,\cdots$. Hence, we have
\[
x_1=\lim_{k\to\infty}\phi_{1,j_1}\circ \phi_{2,j_2}\circ\cdots\circ\phi_{k,j_k}(x_{k+1}).
\]
Using the assumption \eqref{eq1.1}, similar to the proof of Claim 1, we have
\[
x_1=\lim_{k\to\infty}\phi_{1,j_1}\circ \phi_{2,j_2}\circ\cdots\circ\phi_{k,j_k}(x_{k+1})=\lim_{k\to\infty}\phi_{1,j_1}\circ \phi_{2,j_2}\circ\cdots\circ\phi_{k,j_k}(a).
\]
Therefore, $U_1\subseteq K_1$.

\medskip

On the other hand, however, for any given sequence $j_1j_2\cdots \in \Sigma_1^{\mathbb{N}}$, \eqref{eq2.2} implies that there exist points $y_{k+1}\in U_{k+1}$ ($k\ge 1$) such that
\[
\phi_{1,j_1}\circ \phi_{2,j_2}\circ\cdots\circ\phi_{k,j_k}(y_{k+1})\in U_1.
\]
Since $\phi_{n,j}$ are contractive maps on compact set $X$ satisfying \eqref{eq1.1}, we have
\[\begin{array}{rl}
&\rho(\phi_{1,j_1}\circ \phi_{2,j_2}\circ\cdots\circ\phi_{k,j_k}(y_{k+1}),\ \phi_{1,j_1}\circ \phi_{2,j_2}\circ\cdots\circ\phi_{k,j_k}(a))\\
\le& |\phi_{1,j_1}\circ \phi_{2,j_2}\circ\cdots\circ\phi_{k,j_k}(X)|_\rho\to 0
\end{array}\]
as $k\to\infty$. Hence, we have
\[
\lim_{k\to\infty}\phi_{1,j_1}\circ \phi_{2,j_2}\circ\cdots\circ\phi_{k,j_k}(a)=\lim_{k\to\infty}\phi_{1,j_1}\circ \phi_{2,j_2}\circ\cdots\circ\phi_{k,j_k}(y_{k+1})\in U_1.
\]
Therefore $K_1\subseteq U_1$ by using \eqref{eq2.1}. Claim 3 is proven.

\medskip

It is easy to see that all conclusions in (i) follow from Claim 1 $\sim$ 3.

\medskip

(ii)
Consider the symbolic space $\Sigma_n^{\mathbb N}$.  If we define $$\varrho_n(i_ni_{n+1}i_{n+2}\cdots i_{k}\cdots,\ j_nj_{n+1}j_{n+2}\cdots j_{k}\cdots)=0.5^{\min\{k\ge n:\ i_k\ne j_k\}-n},$$ then the pair $(\Sigma_n^{\mathbb N},\ \varrho_n)$ is a compact metric space for each $n>0$.  The topology of $(\Sigma_n^{\mathbb N},\ \varrho_n)$ is also generated by the class of cylinder sets
\[
\left\{\Sigma_n^{\mathbb N}(j_nj_{n+1}j_{n+2}\cdots j_{k}):\ 1\le j_{n+i}\le N_{n+i},\ i\ge 0\right\}.
\]
 Hence, there is a unique Borel probability measure $\nu_n$ on $\Sigma_n^{\mathbb N}$ satisfying
\begin{equation}\label{eq2.5}
\nu_n(\Sigma^{\mathbb{N}}_n(j_nj_{n+1}j_{n+2}\cdots j_{k}))=p_{n,j_n}p_{n+1,j_{n+1}}p_{n+2,j_{n+2}}\cdots p_{k,j_{k}},\quad k\ge n>0.
\end{equation}

\medskip

Note that \eqref{eq2.1} defines a surjective map $\pi_n:\ \Sigma_n^{\mathbb N}\to K_n$ for any fixed $a\in X$. By \eqref{eq1.1} and the compactness of $X$, we see that this $\pi_n$ is a continuous map from $(\Sigma_n^{\mathbb N},\ \varrho_n)$ to $(K_n,\ \rho)$ for each $n>0$. Hence, $\mu_n:=\nu_n\circ\pi_n^{-1}$ is a Borel probability measure on $K_n$ for all $n>0$. Furthermore, \eqref{eq2.5} shows that $\mu_n$ satisfy \eqref{eq2.3} and the support of $\mu_n$ is $K_n$.

\medskip

We then need only to prove the uniqueness of $\mu_n$. Assume that $\left\{\mu'_n\right\}_{n\geq1}$ is another sequence of probability measures satisfying \eqref{eq2.3} and the sequence of their supports is $\left\{K_n\right\}_{n\geq1}$.

\medskip

For any open set $O$ in the metric space $(X,\ \rho)$, by the continuity of $\pi_n$, we see that $\pi_n^{-1}(O)$ is an open set in the metric space $(\Sigma_n^{\mathbb N},\ \varrho_n)$. Hence, $\pi_n^{-1}(O)$ can be written as a disjoint union of some cylinder sets: $\Sigma_n^{\mathbb N}(j_nj_{n+1}j_{n+2}\cdots j_{k})$, i.e.,
\begin{equation}\label{eq2.6}
\pi_n^{-1}(O)=\bigcup_{i}\Sigma_n^{\mathbb N}(J_i),
\end{equation}
where the cylinder sets $\{\Sigma_n^{\mathbb N}(J_i)\}_{i\geq1}$ are disjoint. \eqref{eq2.3} implies $\mu_n'(O) =\sum\limits_{j=1}^{N_n} p_{n,j}\mu_{n+1}' \circ\phi_{n,j}^{-1}(O)$ for all $ n\ge1$. Noting that $O$ is open, \eqref{eq2.6} shows that $\mu_{n+1}' \circ\phi_{n,j}^{-1}(O)>0$ if and only if $j$ is the first symbol of some $J_i$. Hence, by iterating \eqref{eq2.3}, one see that \eqref{eq2.6} shows
\[
\mu_n'(O)=\sum_{i}p_{n,J_i}\mu'_{n+|J_i|}(K_{n+|J_i|})=\sum_{i}p_{n,J_i}=\mu_n(O).
\]
This means that $\mu_n'$ and $\mu_n$ coincide on any open set of $(K,\ \rho)$. Hence, the uniqueness theorem of measure extension implies $\mu'_n=\mu_n$. The uniqueness of $\mu_n$ is proven.
\end{pf}

The above theorem shows that we can give the following definitions.

\begin{define} \label{def2.2} Let $\left\{\Phi_n\right\}_{n=1}^\infty$ with $\Phi_n=\left\{\phi_{n,j}\right\}_{j=1}^{N_n}$ be an MIFS on $(X,\ \rho)$ as in Definition \ref{def1.1}.

\medskip

{\rm (ia)}  If $\phi_{n,j}$ are all similar maps, then $\left\{\Phi_n\right\}_{n=1}^\infty$ is called a Moran-type self-similar iterated function system (Moran-type self-similar IFS or self-similar MIFS).

\medskip
{{\rm (ib)}  If all $\phi_{n,j}$ can be extended to be $C^1$ conformal maps on $X_1$, where $X$ is a compact subset of the interior $X_1^o$, and \eqref{eq1.1} holds if  $\rho(\cdot,\cdot)$ and $X$ are replaced by the Euclidean metric and  $X_1$, respectively,  then $\left\{\Phi_n\right\}_{n=1}^\infty$ is called a Moran-type self-conformal iterated function system (Moran-type self-conformal IFS or self-conformal MIFS).}

\medskip

{\rm (ic)}  If $\phi_{n,j}$ can be extended to an affine map with form $\phi_{n,j}(x)=A^{-1}_{n,j}(x+\alpha_{n,j})$ on the full space ${\mathbb R}^d$ for all $n,\ j$, it is called a Moran-type generalized self-affine iterated function system (Moran-type generalized self-affine IFS or generalized self-affine MIFS).

\medskip

{\rm (id)}  If $A_{n,j}$ in {\rm (ic)} are independent of $j$ for all $n>0$, it is called a Moran-type self-affine iterated function system (Moran-type self-affine IFS or self-affine MIFS).

\medskip

{\rm (ii)}  We call each compact set $K_n$ a Moran-type invariant set. According to different types of MIFSs defined in {\rm (ia)}$\sim${\rm (id)}, we call the $K_n$ a Moran-type self-similar set, a Moran-type self-conformal set, a Moran-type generalized self-affine set and a Moran-type self-affine set, respectively.

\medskip

{\rm (iii)}  We call the probability measure $\mu_n$ a Moran-type invariant measure. According to different types of MIFSs defined in {\rm (ia)}$\sim${\rm (id)}, we call the probability measure $\mu_n$ a Moran-type self-similar measure, a Moran-type self-conformal measure, a Moran-type generalized self-affine measure and a Moran-type self-affine measure, respectively.
\end{define}

\begin{remark}\label{remark2.3} The assumption \eqref{eq1.1} can not be replaced by the assumption that all maps $\phi_{n,j}$ are contractive maps on $(X,\ \rho)$. For example, let $$\Phi_n=\left\{\phi_{n,1}(x)=\left(\frac12\right)^{2^{-n}}x,\ \phi_{n,2}(x)=\frac12(x+1)\right\}$$ on the one dimensional Euclidean space. Then we can choose a compact set $X=[0,\ 1]\subset{\mathbb R}$ such that $\phi_{n,j}(X)\subset X$ for all $n,\ j\geq1$. In this case, the constants $c_{2,n}=(\frac12)^{2^{-n}}<1$ does not satisfy the assumption $\lim\limits_{n\to+\infty}\prod_{j=1}^nc_{2,j}=0$ in \eqref{eq1.1}, but all other assumptions in Definition \ref{def1.1} are satisfied. It is easy to see that
\[
 \lim\limits_{k\to+\infty}\phi_{n,1}\circ\phi_{n+1,1}\cdots\circ\phi_{n+k,1}(a)
\]
is not independent of $a\in X$.
\end{remark}

\begin{remark}\label{remark2.4} The assumption that $(X,\ \rho)$ is a compact metric space can not be replaced by the assumption that $(X,\ \rho)$ is a complete and separable metric space. For example, let $\phi_{n,1}(x)=\frac13 x,\ \phi_{n,2}(x)=\frac13(x+1),\ \phi_{n,3}(x)=\frac13(x+2)$ on the one dimensional Euclidean space $\mathbb R$. Then, the $K_n$ in Definition \ref{def1.1} is the unit interval $[0,\ 1]$ for all $n>0$. But $U_n:=[3^{n-1}a,\ 3^{n-1}a+1]$ for all $a\in {\mathbb R}$ also satisfy
\[
U_{n}=\bigcup_{j=1}^3\phi_{n,j}(U_{n+1}),\quad \phi_{n,i}(U^o_{n+1})\cap \phi_{n,j}(U^o_{n+1})=\emptyset\ (1\le i<j\le 3),\ n>0.
\]
This means that the uniqueness of $K_n$ is invalid if $X$ is not compact.
\end{remark}

\section{Definitions of separation conditions}
\setcounter{equation}{0}\setcounter{theorem}{0}

It is well known that, for self-similar IFS, the open set condition (OSC), the strong separation condition (SSC) and the weak separation condition (WSC) are crucial in the investigations on the self-similar sets or self-similar measures. This section is devoted to extend OSC, SSC and WSC for IFSs to  MIFSs.

\medskip

In order to define separation conditions of MIFSs, we use the idea of  \cite{[DeW]} for generalized self-affine IFSs and the idea of \cite{[LNW]} for self-conformal IFSs.

\medskip
\begin{define} \label{def3.1}  Let $\left\{\Phi_n\right\}_{n=1}^\infty$ with $\Phi_n=\left\{\phi_{n,j}\right\}_{j=1}^{N_n}$ be an MIFS on $(X,\ \rho)$ as in Definition \ref{def1.1}.

\medskip

{\rm (i)}  We say that the MIFS satisfies the Moran-type open set condition (MOSC in short) if there exists a sequence of  open sets $V_n\subset X$ such that
\begin{equation}\label{eq3.1}
\phi_{n,j}(V_{n+1})\subseteq V_{n},\ \ \phi_{n,i}(V_{n+1})\cap \phi_{n,j}(V_{n+1})=\emptyset\ (1\le i\neq j\le N_n),\ \ n\ge1,
\end{equation}
and
\begin{equation}\label{eq3.2}
\inf\left\{{\mathcal L}^d(V_n):\ n=1,2,\cdots\right\}>0.
\end{equation}

\medskip

{\rm (ii)}  We say that the MIFS satisfies the Moran-type weak separation condition (MWSC in short) if there exists a sequence of subsets $U_n\subseteq X$ such that
\begin{equation}\label{eq3.3}
\inf\left\{{\mathcal L}^d(U_n):\ n=1,2,\cdots\right\}>0,
\end{equation}
and
\begin{equation}\label{eq3.4}
\gamma_1:=\sup\limits_{1>b>0}\max\limits_{I\in{\mathcal I}_b}\#\left\{\phi_{1,J}:\ J\in{\mathcal I}_b,\ \phi_{1,J}(U_{|J|+1})\cap\phi_{1,I}(U_{|I|+1})\ne\emptyset\right\}<+\infty.
\end{equation}

{\rm (iii)}  We say that the MIFS satisfies the Moran-type strong separation condition (MSSC in short) if
\[
\phi_{1,J}(K_{n+1})\cap\phi_{1,I}(K_{n+1})=\emptyset
\]
for any $n>0$ and $I\ne J\in{}\Sigma_1^{n}$.
\end{define}

\begin{remark} The ``Moran structure condition" in \cite{[7]} and other references require that all $V_n$ in the definition of MOSC are the same, so MOSC is weaker than the ``Moran structure condition" in  \cite{[7]} and other references. Hence, our Theorem \ref{th4.3} and Theorem \ref{th4.5} in Section 4 are nontrivial generalizations of conclusions in  \cite{[7]} and other references. Some other remarks and examples are given in Section 5 to explain why we introduce the above definitions.
\end{remark}

The idea of \cite[(2.3)]{[DeW]} for generalized self-affine IFSs and the BDP defined in \cite{[LNW]} for self-conformal IFSs are important for investigating attractors. We extend these ideas to introduce the so-called MWHP and MBDP for MIFSs in the following.

\begin{define} \label{def3.3} Let $\left\{\Phi_n\right\}_{n=1}^\infty$ be an MIFS as defined in Definition \ref{def1.1}.
If
\begin{equation}\label{eq3.5}
\gamma_2:=\sup\limits_{1>b>0}\sup\left\{\frac{\rho(\phi_{1,J}(x),\ \phi_{1,J}(y))}{\rho(\phi_{1,I }(x),\ \phi_{1,I}(y))}:\ I,\ J\in {\mathcal I}_b,\ x\ne y\in X\right\}<+\infty,
\end{equation}
we say that the MIFS $\left\{\Phi_n\right\}_{n=1}^\infty$ satisfies the Moran-type weak homogenous property (MWHP in short).

\medskip

If
\begin{equation}\label{eq3.6}
\gamma_3:=\sup\limits_{I\in\Sigma_1^*}\left\{\frac{R_I}{r_I}\right\}<+\infty,
\end{equation}
we say that the MIFS $\left\{\Phi_n\right\}_{n=1}^\infty$ satisfies the Moran-type bounded distortion property (MBDP in short).
\end{define}

\medskip

{\bf Remark}.
If the metric $\rho$ in Definition \ref{def1.1} is translation-invariant, it is easy to see that both self-affine MIFSs and   self-similar MIFSs satisfy the MWHP.

\medskip

For the relationship between MOSC and MWSC, similar to the relationship between OSC (see \cite{[2]}) and WSC (see \cite{[14]}) for IFSs, we have the following conclusion.

\begin{theorem}\label{th3.4} The MOSC implies the MWSC.
\end{theorem}
\begin{pf}
 It is clear that  both \eqref{eq3.3} and \eqref{eq3.4} hold when we take $U_n$ to be the $V_n$ in \eqref{eq3.1} and \eqref{eq3.2}. Hence, the MOSC implies the MWSC.

\end{pf}

  By the fact that the MOSC ensures that $\phi_{1,I}\neq \phi_{1,J}$ when $I\neq J$, we have the following lemma which plays a key role in the dimension theory of Moran-type self-similar sets as shown in Section 4.

\begin{lemma}\label{th3.5}
Assume $\left\{\Phi_n\right\}_{n=1}^\infty$ is a self-conformal MIFS satisfying the MBDP and MWHP or a generalized self-affine MIFS satisfying the MWHP.

\medskip

 {\rm (i)} If $\left\{\Phi_n\right\}_{n=1}^\infty$ satisfies the MWSC, then
\begin{equation}\label{eq3.7}
\gamma_4:=\sup\limits_{1>b>0}\max\limits_{I\in{\mathcal I}_b}\#\left\{\phi_{1,J}\in{\mathcal A}_b:\ \phi_{1,J}(X)\cap\phi_{1,I}(X)\ne\emptyset\right\}<+\infty.
\end{equation}

\medskip

{\rm (ii)} If $\left\{\Phi_n\right\}_{n=1}^\infty$ satisfies the MOSC, then
\begin{equation}\label{eq3.8}
\gamma'_4:=\sup\limits_{1>b>0}\max\limits_{I\in{\mathcal I}_b}\#\left\{J\in{\mathcal I}_b:\ \phi_{1,J}(X)\cap\phi_{1,I}(X)\ne\emptyset\right\}<+\infty.
\end{equation}
\end{lemma}

{\bf Remark}. Example \ref{Ex5.6} shows that, without the MWHP, the conclusion of Lemma \ref{th3.5} maybe wrong. This means that MWSC for MIFSs is much different from WSC for IFSs.

\medskip
\begin{pf} (i)
For the case that $\left\{\Phi_n\right\}_{n=1}^\infty$ is a generalized self-affine MIFS satisfying the MWHP, every $\phi_{n,j}$ is well defined on ${\mathbb R}^d$. For any $I\in{\mathcal I}_b$, let $\left\{\phi_{1,J}\in{\mathcal A}_b:\ \phi_{1,J}(X)\right.$ $\left.\cap\phi_{1,I}(X)\ne\emptyset\right\}=\{\phi_{1,J_1},\ \phi_{1,J_2},\ \cdots,\ \phi_{1,J_{n_b}}\}$. Then the MWSC shows there exists a sub family $\{\phi_{1,J_{i_1}},$ $\phi_{1,J_{i_2}},$ $\cdots,\ \phi_{1,J_{i_k}}\}\subseteq\{\phi_{1,J_1},\ \phi_{1,J_2},\ \cdots,\ \phi_{1,J_{n_b}}\}$ such that $k\ge n_b\gamma_1^{-1}$ and $\{\phi_{1,J_{i_1}}(U_{|J_{i_1}|+1}),$ $\phi_{1,J_{i_2}}(U_{|J_{i_2}|+1}),$ $\cdots,$ $\phi_{1,J_{i_k}}(U_{|J_{i_k}|+1})\}$ are disjoint.
The MWHP shows
\[
\max\{\|A_{1,J_{i_j}}A^{-1}_{1,I}\|,\ \|A_{1,I}A_{1,J_{i_j}}^{-1}\|:\ j=1,2,\cdots,k\}\le \gamma_2.
\]
Hence, we have $X_{\gamma_2|X|}\supset\bigcup\limits_{j=1}^k\phi_{1,I}^{-1}\circ\phi_{1,J_{i_j}}(X)\supset\bigcup\limits_{j=1}^k\phi_{1,I}^{-1}\circ\phi_{1,J_{i_j}}(U_{|J_{i_j}|+1})$, where $X_{\delta}=\{x\in{\mathbb R}^d:\ \mbox{dist}(x,\ X)<\delta\}$ is the $\delta$-neighborhood of $X$. This means
\[\begin{array}{rl}
{\mathcal L}^d(X_{\gamma_2|X|})\ge&\sum_{j=1}^k{\mathcal L}^d(\phi_{1,I}^{-1}\circ\phi_{1,J_{i_j}}(U_{|J_{i_j}|+1}))\\ \ge& \gamma_2^{-d}\sum_{j=1}^k{\mathcal L}^d(U_{|J_{i_j}|+1})\\
\geq&k\gamma_2^{-d}\inf\{{\mathcal L}^d(U_{n}):\ n>0\} .
\end{array}
\]
Since the MWSC guarantees that $\inf\{{\mathcal L}^d(U_{n}):\ n>0\}>0$, we get $k\le \frac{\gamma_2^d{\mathcal L}^d(X_{\gamma_2|X|})}{\inf\{{\mathcal L}^d(U_{n}):\ n>0\}}$. It follows that $n_b\le \frac{\gamma_1\gamma_2^d{\mathcal L}^d(X_{\gamma_2|X|})}{\inf\{{\mathcal L}^d(U_{n}):\ n>0\}}$, which implies that the upper-bound of $n_b$ for all $b\in(0,\ 1)$ and  $I\in{\mathcal I}_b$ is finite.  Therefore, \eqref{eq3.7} holds.

\medskip

For the case that $\left\{\Phi_n\right\}_{n=1}^\infty$ is a self-conformal MIFS satisfying the MBDP and MWHP, the proof is similar and  we omit the details.

\medskip

(ii) The proof is similar to (i) if we replace $\phi_{1,J}$ by $J$.
\end{pf}

\medskip

The following theorem shows that the above MWHP is really an extension of \cite[(2.3)]{[DeW]} for generalized self-affine IFSs and the above MBDP is really an extension of the BDP defined in \cite{[LNW]} for self-conformal IFSs.

\medskip

\begin{theorem}\label{th3.6} Let $\Phi=\{\phi_j\}_{j=1}^m$ be an IFS  and $\Phi_{n}=\Phi$ for all $n\ge1$ on the compact metric space $(X,\ \rho)$.

\medskip

{\rm (i)} If $\Phi$ is a generalized self-affine IFS, then $\Phi$ satisfies the \cite[(2.3)]{[DeW]} if and only if $\left\{\Phi_n\right\}_{n=1}^\infty$ satisfies the MWHP.

\medskip

{\rm (ii)} If $\Phi$ is a self-conformal IFS, then $\Phi$ satisfies the BDP defined by \cite[(1.7)]{[LNW]} if and only if $\left\{\Phi_n\right\}_{n=1}^\infty$ satisfies the MBDP.
\end{theorem}

\begin{pf} The assumption that $\Phi_{n}=\Phi$ for all $n\ge1$ shows that $N_n=m$ and $\phi_{n,j}=\phi_j$ for all $n\ge1$ and $N_n\ge j\geq1$.

\medskip

(i)  Assume that all maps in $\Phi$ have forms
\[
\phi_j(x)=A^{-1}_j(x+\alpha_j),\quad j=1,2,\cdots,m,\quad x\in\mathbb{R}^d,
\]
where $\alpha_j \in\mathbb{R}^d$ and $A_j$ are expansive matrices satisfying $c_1\rho(x,\ y)\le\rho(A^{-1}_jx,\ A^{-1}_jy)\le c_2\rho(x,\ y)$ with $0<c_1\le c_2<1$.
We write $A_{j_1 j_2\cdots j_n}=A_{j_n}\cdots A_{ j_2}A_{j_1}$ (so $A^{-1}_{j_1 j_2\cdots j_n}=A^{-1}_{j_1}A^{-1}_{ j_2}\cdots A^{-1}_{j_n}$), $\rho_j=|\det(A_j)|^{-1/d}$, $\underline{\rho}=\min\{\rho_j:\ 1\le j\le m\}$ and $\rho_{j_1 j_2\cdots j_n}=|\det(A^{-1}_{j_1}A^{-1}_{ j_2}\cdots A^{-1}_{j_n})|^{1/d}$.

\medskip

It is easy to see that, without loss of generality, one can replace both the $X$ in the definition of MIFS $\left\{\Phi_n\right\}_{n=1}^\infty$ and the $K_{\varepsilon_0}$ in the definition of $\Gamma_b$ on page 1355 of \cite{[DeW]} by a sufficiently large ball $B(a,\ r)$.
Hence, we need only to prove that \eqref{eq3.5} holds for the ${\mathcal I}_b$ defined in \eqref{eq1.4} if and only if it holds when ${\mathcal I}_b$ is replaced by the $\Gamma_b$ defined on page 1355 of \cite{[DeW]}.

\medskip

Also, we have
\[
\frac{\|\phi_{J}(x)-\phi_{J}(y)\|}{\|\phi_{I }(x)-\phi_{I}(y)\|}=\frac{\|A^{-1}_{J}(x-y)\|}{\|A^{-1}_{I }(x-y)\|}=\frac{\|A^{-1}_{J}A_{I }\alpha\|}{\|\alpha\|},
\]
where $\alpha=A^{-1}_{I }(x-y)$ and $\|\cdot\|$ is the Euclidean norm. Note that, for the Euclidean norm $\|\cdot\|$ and the  metric $\rho$, there is a constant $\alpha_0\ge 1$ such that
\begin{equation}\label{eq3.9-1}
\alpha_0^{-1}\le\frac{\|x-y\|}{\rho(x,\ y)}\leq \alpha_0,\quad\forall\ x\ne y\in X.
\end{equation}
Hence, \eqref{eq3.5} is equivalent to
\begin{equation}\label{eq3.10}
L:=\sup\limits_{1>b>0}\sup\left\{\|A^{-1}_{J}A_{I }\|:\ I,\ J\in {\mathcal I}_b\right\}<+\infty.
\end{equation}
Recall that the $\Gamma_b$ is  defined on page 1355 of \cite{[DeW]} by
\[
\Gamma_b=\left\{j_1 j_2 \cdots j_n\in{}\Sigma_1^*:\ \rho_{j_1 j_2\cdots j_n} \le b < \rho_{j_1 j_2\cdots j_{n-1}}\right\}.
\]
The expression (2.3) in \cite{[DeW]} is
\begin{equation}\label{eq3.11}
\tau:=\sup\{\|A_IA^{-1}_J\|:\ I,\ J\in\Gamma_b,\ 0<b<1\}<+\infty.
\end{equation}
In the following, we will prove that \eqref{eq3.10} and \eqref{eq3.11} are equivalent.

\medskip

Assume that \eqref{eq3.10} holds with the ${\mathcal I}_b$ defined in \eqref{eq1.4}. Fix  $b_0\in(0,\ 1)$, and  $I,\ J\in\Gamma_{b_0}$.
Then, we obtain  $I\in{\mathcal I}_{a_0}$ and $J\in{\mathcal I}_{a_1}$ for some $a_0,\ a_1\in(0,\ 1)$. Without loss of generality, suppose $a_0\ge a_1$. There exist finite words $J_1\in{\mathcal I}_{a_0}$ and $J_0\in{}\Sigma^*_{1}$ such that $J=J_0 J_1$. Hence, \eqref{eq3.10} shows
\begin{equation}\label{eq3.12}
\|A^{-1}_{I}A_{J_1}\|\le L_0,\quad  \|A^{-1}_{J_1}A_{I}\|\le L_0
\end{equation}
for some constant $L_0>0$.
As a result, all singular values of $A^{-1}_{I}A_{J_1}$ belong to the interval $[L_0^{-1},\ L_0]$, so $L_0^{-d}\le|\det(A^{-1}_{I})\det(A_{J_1})|\le L_0^d$. Consequently, we have
\[
(\underline{\rho})^dL_0^{-d}{b_0}^d\le L_0^{-d}|\det(A^{-1}_{I})|\le |\det(A^{-1}_{J_1})|\le L_0^d|\det(A^{-1}_{I})|\le L_0^d{b_0}^d
\]
by using $I\in\Gamma_{b_0}$. Since $J\in\Gamma_{b_0}$, the definition of $\Gamma_{b_0}$ shows $|\det(A^{-1}_{J_1})\det(A^{-1}_{J_0})|$ $=|\det(A^{-1}_{J})|$ belongs to $[(\underline{\rho})^d{b_0}^d,\ {b_0}^d]$. Hence, $|\det(A^{-1}_{J_0})|\ge (\underline{\rho})^dL_0^{-d}$, which implies  that the length of $J_0$ is bounded by a constant independent of $b_0$ and $I,\ J$. Therefore, the number of  all such $A_{J_0}$ for all $b_0$ and $I,\ J$ is finite. Noting that all $A_j$ are invertible matrices, we see that there is a constant $c\ge1$ such that
$c^{-1}\le \|A_{J_0}\|,\ \|A^{-1}_{J_0}\|\le c$. Thus, \eqref{eq3.12} implies
\[
\|A^{-1}_{I}A_{J}\|\le cL_0,\quad  \|A^{-1}_{J}A_{I}\|\le cL_0.
\]
Since $cL_0$ is a constant, one see that \eqref{eq3.11} holds.

\medskip

Now we assume that  \eqref{eq3.11} holds.
\medskip
 By \eqref{eq3.9-1}, for any $I\in \Sigma_1^*$, we have
 \begin{equation}\label{metric-norm}
  \frac{1}{\alpha_0^2}\|A_I^{-1}\|\leq R_I\leq \alpha_0^2\|A_I^{-1}\|.
 \end{equation}
Fix a number $a_0\in(0,\ 1)$, and  $I,\ J\in{\mathcal I}_{a_0}$. According to  the definition of ${\mathcal I}_{a_0}$ and \eqref{metric-norm}, we have
\begin{equation}\label{eq3.13-1}
\alpha_0^2a_0\geq \|A_I^{-1}\|\geq \frac{1}{\alpha_0^2}R_I\geq \frac{1}{\alpha_0^2}R_{I^-}\min_{1\le j\le m}\{R_j\}\geq \frac{a_0}{\alpha_0^4}\min_{1\le j\le m}\{\|A_j^{-1}\|\},
\end{equation}
where $I^-$ denotes the finite word obtained by deleting the last digit of $I$.
Without loss of generality, we assume that  $b_1\ge b_2$. There exist $J_1\in\Gamma_{b_1}$ and $J_2\in{}\Sigma_{|J_1|+1}^*$ such that $J=J_1J_2$. 
 By means of \eqref{eq3.11}, we get $\|A_IA^{-1}_{J_1}\|\le\tau$ and  $\|A_{J_1}A^{-1}_I\|\le\tau$. Hence, $\tau^{-1}\|A^{-1}_I\|\le\|A^{-1}_{J_1}\|\le\tau\|A^{-1}_I\|$. Combining \eqref{metric-norm} and \eqref{eq3.13-1} shows
 $$\|A^{-1}_{J_2}\|\ge \|A^{-1}_{J_1}\|^{-1}\|A^{-1}_{J}\|\ge   \frac{1}{\tau\|A_I^{-1}\|}\frac{a_0}{\alpha^4_0}\min_{1\le j\le m}\{\|A_j^{-1}\|\}\ge \frac{1}{\tau\alpha_0^6}\min_{1\le j\le m}\{\|A_j^{-1}\|\}.$$
 This means that the length of $J_2$ is bounded by a constant independent of $a_0$ and $I,\ J\in{\mathcal I}_{a_0}$. Therefore, the number of all such $A_{J_2}$ for all $a_0$ and $I,\ J\in{\mathcal I}_{a_0}$ is finite. Noting that all $A_j$ are invertible matrices, we see that there is a constant $c\ge1$ such that
$c^{-1}\le \|A_{J_2}\|,\ \|A^{-1}_{J_2}\|\le c$.
 Then, combining $\|A_IA^{-1}_{J_1}\|\le\tau$ and  $\|A_{J_1}A^{-1}_I\|\le\tau$ shows that
 \eqref{eq3.10} holds. Conclusion (i) is proven.

\medskip

(ii) When $\Phi$ is a self-conformal IFS,  we have $\rho(x,\ y)=\|x-y\|$. Since the set $V$ in \cite[(1.7)]{[LNW]} is actually the $X_1^o$, we see that the definitions of ${\mathcal I}_b$ and ${\mathcal A}_b$ coincide with \eqref{eq1.4}.

\medskip

Recall that $\sup\{\frac{\|\phi_{J}(x)-\phi_{J}(y)\|}{\|x-y\|}:\ x\ne y\in X\}=R_J$ and $\inf\{\frac{\|\phi_{J}(x)-\phi_{J}(y)\|}{\|x-y\|}:\ x\ne y\in X\}=r_J$ for all $J\in\Sigma^*$. We see that the BDP for $\Phi$ and the MBDP for $\{\Phi_n\}_{n=1}^\infty$ are equivalent.
\end{pf}

The following theorem shows that MOSC, MWSC and MSSC for MIFSs are really extensions of OSC, WSC and SSC both for generalized self-affine IFSs and self-conformal IFSs.

\begin{theorem}\label{th3.7} Let  $\Phi_{n}=\Phi$ for all $n\ge1$. Suppose $\Phi$ is a generalized self-affine IFS satisfying \cite[(2.3)]{[DeW]} or a self-conformal IFS satisfying the BDP. Then, $K_n=K_1$ for all $n\ge1$ and the following three statements hold.

\medskip

\rm{(i)} The MIFS $\left\{\Phi_n\right\}_{n=1}^\infty$ satisfies the MOSC if and only if the IFS $\Phi$ satisfies the  OSC.

\medskip

\rm{(ii)} The MIFS $\left\{\Phi_n\right\}_{n=1}^\infty$ satisfies the MWSC if and only if the IFS $\Phi$ satisfies the WSC.

\medskip

\rm{(iii)} The MIFS $\left\{\Phi_n\right\}_{n=1}^\infty$ satisfies the MSSC if and only if the IFS $\Phi$ satisfies the  SSC.
\end{theorem}
\begin{pf} It is obvious that $K_n=K_1$ for all $n\ge1$. We need only to prove conclusions (i) $\sim$ (iii).

\medskip

(i) By taking $V_n=V$ in Definition \ref{def3.1}, the  necessity is obvious.

\medskip

Suppose $\left\{\Phi_n\right\}_{n=1}^\infty$ satisfies the MOSC. Then, there exists a sequence of open sets $V_n$ such that \eqref{eq3.1} and \eqref{eq3.2} hold. Since $\Phi_n=\Phi$ for all $n>0$, it is easy to see that the MOSC of $\left\{\Phi_n\right\}_{n=1}^\infty$  implies that $\phi_J\ne \phi_J$ when $I\ne J$.

\medskip

We divide the proof of the sufficiency into two parts.

\medskip

{\bf Part 1}. $\Phi$ is a generalized self-affine IFS satisfying the \cite[(2.3)]{[DeW]}.

\medskip

Let $\phi_j(x)=A_j^{-1}(x+\alpha_j)$, $j=1,2,\dotsc, m$. Thus $\phi_j$ are well defined in ${\mathbb R}^d$. Recall
\[
\Gamma_b=\{j_1j_2\cdots j_s\in{}\Sigma_1^*:\ |\det(A_{j_1j_2\cdots j_s}^{-1})|\le b^d<|\det(A_{j_1j_2\cdots j_{s-1}}^{-1})|\},\quad b\in(0,\ 1).
\]
We first prove that the IFS $\Phi$ satisfies
\begin{equation}\label{eq3.14}
\max_{0<b<1}\max_{I\in\Gamma_b}\#\left\{J\in\Gamma_b:\ \phi_{J}(X)\cap\phi_{I}(X)\ne\emptyset\right\}<+\infty.
\end{equation}
For any $b\in(0,\ 1)$ and $I\in\Gamma_b$, let $\{J\in\Gamma_b:\ \phi_J(X)\cap\phi_I(X)\ne\emptyset\}=\{J_1,\ J_2,\ \dotsc,\ J_q\}$. Then, we obtain
\[
X\cap\phi_I^{-1}\circ\phi_{J_j}(X)\ne\emptyset,\quad j=1,\ 2,\ \cdots,\ q.
\]
According to  the assumption \cite[(2.3)]{[DeW]}, there is a constant $L$ independent of $b\in(0,\ 1)$ and $I\in\Gamma_b$ such that
$L^{-1}\|x-y\|\le\|A_{J_j}^{-1}A_I(x-y)\|\le L\|x-y\|$.
Hence, we have
\begin{equation}\label{eq3.15}
{\mathcal L}^d(\phi_I^{-1}\circ\phi_{J_j}(V_{|J_j|+1}))\ge L^{-d}\inf_{n\ge1}\{{\mathcal L}^d(V_n)\},
\end{equation}
\begin{equation}\label{eq3.16}
\phi_I^{-1}\circ\phi_{J_j}(V_{|J_j|+1})\subset\phi_I^{-1}\circ\phi_{J_j}(X)\subset B(x_0, (L+1)|X|)\quad j=1,\ 2,\ \cdots,\ q,
\end{equation}
where $x_0$ is an arbitrary point of $X$. The MOSC shows that $\phi_{J_j}(V_{|J_j|+1})$ are disjoint, so $\phi_I^{-1}\circ\phi_{J_j}(V_{|J_j|+1})$ are disjoint. Hence, \eqref{eq3.15} and \eqref{eq3.16} imply that
\[
qL^{-d}\inf_{n\ge1}\{{\mathcal L}^d(V_n)\}\le \sum_{j=1}^q{\mathcal L}^d(\phi_I^{-1}\circ\phi_{J_j}(V_{|J_j|+1}))\le{\mathcal L}^d(B(x_0, (1+L)|X|)),
\]
Therefore, the assumption \eqref{eq3.2} in the definition of MOSC shows that $q$ is bounded uniformly, so \eqref{eq3.14} holds.

\medskip

Now, we turn to the OSC. Let $U=X^o$ be the interior of $X$. \eqref{eq3.14} shows
\[
p:=\max_{b>0}\max_{I\in\Gamma_b}\#\left\{{J}:\ J\in\Gamma_b,\ \phi_{J}(U)\cap\phi_{I}(U)\ne\emptyset\right\}<+\infty.
\]
{Hence, there exists a positive number
$b>0$ and $I\in\Gamma_b$ such that 
\[\left\{{J}:\ J\in\Gamma_b,\ \phi_{J}(U)\cap\phi_{I}(U)\ne\emptyset\right\}=\left\{J_1,J_2,\cdots,J_p\right\}.\]
We define
\[
V=\bigcup_{\sigma\in\Sigma^*}\phi_{\sigma}(\phi_I(U)),
\]
where $\Sigma^*$ denotes the group of finite words on the symbolic space $\{1,2,\cdots,m\}^\mathbb{N}.$
 It is clear that $V$ is an open set and $\phi_j(V)\subseteq V$ for all $1\leq j\leq m$. Assume that $\phi_i(V)\cap\phi_j(V)\ne\emptyset$ for some $i\ne j$, i.e.,  $\phi_1(V)\cap\phi_2(V)\ne\emptyset$ without loss of generality. From the definition  of $V$, it follows that there exist two finite words $\sigma$ and $\tau\in\Sigma^*$ such that
  \begin{equation}\label{non-inters}
  \phi_{1\sigma I}(U)\cap\phi_{2\tau I}(U)\ne\emptyset.
   \end{equation}  Let $a=|\det(A^{-1}_{1\sigma })|^{1/d}$ and $a'=|\det({}A^{-1}_{2\tau })|^{1/d}.$
   The definition of $\Gamma_b$ shows that $1\sigma I\in\Gamma_{ab}$ and $2\tau I\in\Gamma_{a'b}$. 
   Without loss of generality, we assume $a\ge a'$. Then, $2\tau I$ can be decomposed as $2\tau I=2\tau_1\tau_2$ with  $2\tau_1\in{\Gamma_{ab}}$ and $\tau_2\in\Sigma^*$. In combination with the fact  $\phi_j(U)\subseteq U$ for all $1\leq j\leq m$ and \eqref{non-inters}, we have $\phi_{1\sigma I}(U)\cap\phi_{2\tau_1}(U)\ne\emptyset$.
Hence, we have
\[
\left\{{J}:\ J\in\Gamma_{ab},\ \phi_{J}(U)\cap\phi_{1\sigma I}(U)\ne\emptyset\right\}\supseteq \left\{1\sigma J_1,1\sigma J_2,\cdots,1\sigma J_p\right\}\cup\left\{2\tau_1\right\}.
\]
This means $\#\left\{{J}:\ J\in\Gamma_{ab},\ \phi_{J}(U)\cap\phi_{1\sigma I}(U)\ne\emptyset\right\}>p$, a contradiction to the definition of $p$.
Hence, $\Phi$ satisfies the OSC with the open set $V$. The statement (i) is proven when $\Phi$ is a generalized self-affine IFS.}

\medskip

{\bf Part 2}. $\Phi$ is a self-conformal IFS satisfying the BDP.

\medskip

Assume that  the MIFS $\left\{\Phi_n\right\}_{n=1}^\infty$ satisfies the MOSC with open sets $\{V_n\}_{n=1}^\infty$.
For any $b\in(0,\ 1)$ and $I\in {\mathcal I}_b$, let
\[
\{J\in {\mathcal I}_b:\ \phi_{J}(X)\cap\phi_{I}(X)\neq \emptyset\}=\{J_1,\ J_2,\ \cdots,\ J_s\}.
\]
Note $|\phi_{J}(X)|\leq b|X|$ for any $J\in {\mathcal I}_b$. By taking a point $x_0\in X$, we have
\[
\phi_{J_i}(V_{|J_i|+1})\subset\phi_{J_i}(X)\subset B(\phi_{I}(x_0),2b|X|), \quad i=1,\ 2,\ \cdots,\ s.
\]
The MOSC shows that $\phi_{J_i}(V_{|J_i|+1})$ are disjoint. According to the BDP, for any   $J\in {\mathcal I}_b$ we have
\[\begin{array}{rl}
b|x-y|\ge&|\phi_{J}(x)-\phi_{J}(y)|\\ \ge& r_{J}|x-y|\\ \ge&\gamma_3^{-1}R_{J^-}\min\{r_j:\ 1 \le j\le m\}|x-y|\\ \ge&\gamma_3^{-1}b\min\{r_j:\ 1\le j\le m\}|x-y|,
\end{array}
\] where $J^-$ denotes the finite word by deleting the last digit of $J$. Hence, we have
\[
\begin{array}{rl}
s(\gamma_3^{-1}b\min\{r_j:\ 1\le j\le m\})^{d}\inf_{n\ge1}\{{\mathcal L}^d(V_n)\}
\le & \sum_{j=1}^s{\mathcal L}^d(\phi_{J_j}(V_{|J_j|+1}))\\ \le &b^d{\mathcal L}^d(B(\phi_{I}(x_0),2|X|)).
\end{array}
\]
Therefore, the upper bound of  $s$ is independent of  $b\in(0,\ 1)$, which implies \eqref{eq3.14} holds. Hence, the IFS $\Phi$ satisfies the WSC. Recall  that the MOSC shows that $\phi_I\ne\phi_J$ when $I\ne J$. As a consequence of  \cite[Theorem 1.3]{[DN1]}, we get  that the IFS $\Phi$ satisfies the OSC.

\medskip

Conclusion (i) is proven.

\medskip

(ii)  Since $\Phi_n=\Phi$ and we have proven, in the proof of Theorem \ref{th3.6}, that \eqref{eq3.10} and \eqref{eq3.11} are equivalent, we see that the WSC will not change when the ${\mathcal I}_b$ is replaced by $\Gamma_b$ in the definition of WSC. Hence the definitions of the MWSC of MIFS $\{\Phi_n\}$ and the WSC of IFS $\Phi$ in \cite{[DeW]} are the same. Finally, a similar argument as the proof of (i) proves the conclusion (ii).

\medskip

(iii) By $\Phi_n=\Phi$ and $K_n=K_1$ for all $n\geq1$, the definitions of the MSSC for MIFS $\{\Phi_n\}_{n\geq1}$ and the SSC for IFS $\Phi$ are the same.
\end{pf}

Furthermore, similar to the OSC and WSC for IFSs, we have the following characterizations.

\begin{theorem}\label{th3.8}
Let $\left\{\Phi_n\right\}_{n=1}^\infty$ be a self-conformal MIFS satisfying the MWHP and MBDP defined in Definition \ref{def3.3} or generalized self-affine MIFS  satisfying the MWHP.

\medskip

{\rm (i)} The MOSC does not hold if and only if there exists a sequence of numbers $b_n\in (0,\ 1)$, a sequence of words $\sigma_n\in{\mathcal I}_{b_n}$ and a sequence of pairs of words
\begin{equation}\label{eqww8-23}
{\mathcal B}_n\subseteq\left\{(I, J): I\ne J\in{\mathcal I}_{b_n},\ \phi_{1,I}(X)\cap\phi_{1,\sigma_n}(X)\ne\emptyset,\ \phi_{1,J}(X)\cap\phi_{1,\sigma_n}(X)\ne\emptyset\right\}
\end{equation}
such that
\[
\lim_{k\to\infty}\sup_{(I, J)\in{\mathcal B}_k}\sup_{x\in X}\|\phi_{1,I}^{-1}\circ\phi_{1,J}(x)-x\|=0,\quad \lim_{k\to\infty}\#{\mathcal B}_k=+\infty.
\]

\medskip

{\rm (ii)} The MWSC does not hold if and only if there exists a sequence of numbers $b_n\in (0,\ 1)$, a sequence of words $\sigma_n\in{\mathcal I}_{b_n}$ and a sequence of pairs of maps
\[
{\mathcal C}_n\subseteq\left\{(\phi_{1,I}, \phi_{1,J}): \phi_{1,I}\ne \phi_{1,J}\in{\mathcal A}_{b_n}, \phi_{1,I}(X)\cap\phi_{1,\sigma_n}(X)\ne\emptyset, \phi_{1,J}(X)\cap\phi_{1,\sigma_n}(X)\ne\emptyset\right\}
\]
such that
\[
\lim_{k\to\infty}\sup_{(f, g)\in{\mathcal C}_k}\sup_{x\in X}\|f^{-1}\circ g(x)-x\|=0,\quad \lim_{k\to\infty}\#{\mathcal C}_k=+\infty.
\]
\end{theorem}
\begin{pf}
Noting $X\subset {\mathbb R}^d$, it is easy to see that the metric $\rho$ in the definition of ${\mathcal I}_b$ can be replaced by the Euclidean norm  $\|\cdot\|$ in order to prove our theorem. We will not use $\rho$ in the following argument.

\medskip

(i) We first prove the necessity of (i) and  assume that the MOSC does not hold.

\medskip

As a consequence of Lemma \ref{th3.5},  there exists a decreasing sequence of positive numbers $b_n\in(0,1)$ and a sequence of finite words  $\sigma_n\in{\mathcal I}_{b_n}$ such that
\begin{equation}\label{eq3.17}
\lim\limits_{n\to+\infty}\#\left\{J\in{\mathcal I}_{b_n}:\ \phi_{1,J}(X)\cap\phi_{1,\sigma_n}(X)\ne\emptyset\right\}=+\infty.
\end{equation}
Fix $x_0\in X$. For $J\in{\mathcal I}_{b_n}$, we  define
\[
S_J(x)=\left\{\begin{array}{ll}
\phi_{1,\sigma_n}(x_0)+b_n^{-1}(\phi_{1,J}(x)-\phi_{1,\sigma_n}(x_0)),& \{\Phi_n\}_{n\geq1}\mbox{ is self-conformal},\\
\phi_{1,\sigma_n}(x_0)+A_{1,\sigma_n}(\phi_{1,J}(x)-\phi_{1,\sigma_n}(x_0)),& \{\Phi_n\}_{n\geq1}\mbox{ is generalized self-afine}.
\end{array}\right.
\]
 We  write
\[
Y=\bigcup_{n=1}^\infty\bigcup\left\{S_{J}(X):\ J\in{\mathcal I}_{b_n}\right\}.
\]
It is clear that $Y$ is bounded since $X$ is compact and $S_J(X)\subset X$ for any  $J\in \Sigma_1^*.$
Given an integer vector $(j_1,\cdots, j_d)\in \mathbb{Z}^d$, let $E_{j_1,j_2,\cdots,j_d}$ stand for the cube
\[
E_{j_1,j_2,\cdots,j_d}:=\prod_{i=1}^d[j_i/k,\ (j_i+1)/k).
\]
Fix an integer $k>0$. By means of the cubes, we get two coverings:
\[
{\mathcal X}_k=\left\{E_{j_1,j_2,\cdots,j_d}:\ X\cap E_{j_1,j_2,\cdots,j_d}\ne\emptyset,\ j_1,j_2,\cdots,j_d\in{\mathbb Z}\right\},
\]
and
\[
{\mathcal Y}_k=\left\{E_{j_1,j_2,\cdots,j_d}:\ Y\cap E_{j_1,j_2,\cdots,j_d}\ne\emptyset,\ j_1,j_2,\cdots,j_d\in{\mathbb Z}\right\}.
\]
It is easy to see that both  ${\mathcal X}_k$ and $\mathcal{Y}_k$ are finite coverings. Fix a point $x_{j_1,j_2,\cdots,j_d}\in X\cap E_{j_1,j_2,\cdots,j_d}$ for each $E_{j_1,j_2,\cdots,j_d}\in{\mathcal X}_k$.
Thanks to \eqref{eq3.17},  there exists an integer $n_k$ such that
\[
\#\left\{J\in{\mathcal I}_{b_{n_k}}:\ \phi_{1,J}(X)\cap\phi_{1,\sigma_{n_k}}(X)\ne\emptyset\right\}>k[\#{\mathcal Y}_k]^{[\#{\mathcal X}_k]}.
\]
Note a fact that for any $J\in{\mathcal I}_{b_{n_k}}$, the point $S_J(x_{j_1,j_2,\cdots,j_d})$ must belong to a unique cube in $\mathcal{Y}_k$.
Hence, there exist
\[
{\mathcal B}_k\subseteq\left\{(I,\ J): I\ne J\in{\mathcal I}_{b_{n_k}},\ \phi_{1,I}(X)\cap\phi_{1,\sigma_{n_k}}(X)\ne\emptyset,\ \phi_{1,J}(X)\cap\phi_{1,\sigma_{n_k}}(X)\ne\emptyset\right\}
\]
such that (a) $\#{\mathcal B}_k>k$; (b) for any pair $(I, J)\in{\mathcal B}_k$, both $S_{I}(x_{j_1,j_2,\cdots,j_d})$ and $S_{J}(x_{j_1,j_2,\cdots,j_d})$ belong to a common cube in $\mathcal{Y}_k$. 
  Hence, for any pair $(I, J)\in{\mathcal B}_k$, we have $\|S_{I}(x_{j_1,j_2,\cdots,j_d})-S_{J}(x_{j_1,j_2,\cdots,j_d})\|\le \sqrt{d}/k$. For any  $x\in E_{j_1,j_2,\cdots,j_d}$, we see $\|S_{J}(x)-S_{J}(x_{j_1,j_2,\cdots,j_d})\|=b_{n}^{-1}\|\phi_{J}(x)-\phi_{J}(x_{j_1,j_2,\cdots,j_d})\|  \le\|x-x_{j_1,j_2,\cdots,j_d}\|\le \sqrt{d}/k$. Similarly, we have $\|S_{I}(x)-S_{I}(x_{j_1,j_2,\cdots,j_d})\|
\le \sqrt{d}/k$. It follows that for any $x\in X$, we have
\[
\|S_{I}(x)-S_{J}(x)\|=b_{n_k}^{-1}\|\phi_{I}(x)-\phi_{J}(x)\|\le 3\sqrt{d}/k.
\]

The above argument shows that there exist $n_k>0$, $\sigma_{n_k}\in {\mathcal I}_{b_{n_k}}$ and
\[
{\mathcal B}_k\subseteq\left\{(I,\ J): I\ne J\in{\mathcal I}_{b_{n_k}},\ \phi_{1,I}(X)\cap\phi_{1,\sigma_{n_k}}(X)\ne\emptyset,\ \phi_{1,J}(X)\cap\phi_{1,\sigma_{n_k}}(X)\ne\emptyset\right\}
\]
such that
\[
\lim_{k\to\infty}\sup_{(I, J)\in{\mathcal B}_k}\sup_{x\in X}\|S_{I}(x)-S_{J}(x)\|=0,\quad \lim_{k\to\infty}\#{\mathcal B}_k=+\infty,
\]
which implies that \[
\lim_{k\to\infty}\sup_{(I, J)\in{\mathcal B}_k}\sup_{x\in X}b_{n_k}^{-1}\|\phi_{1,I}(x)-\phi_{1,J}(x)\|=0,\quad \lim_{k\to\infty}\#{\mathcal B}_k=+\infty.
\]

If $\{\Phi_n\}_{n\geq1}$ is a self-conformal MIFS, by the MWHP, MBDP and the definition of $S_J$, we see that
$b_{n_k}(\gamma_2\gamma_3)^{-2}\le\frac{\|\phi_{1,J}(x)-\phi_{1,J}(y)\|}{\|x-y\|}\le b_{n_k}(\gamma_2\gamma_3)^2$ and $b_{n_k}(\gamma_2\gamma_3)^{-2}\le\frac{\|\phi_{1,I}(x)-\phi_{1,I}(y)\|}{\|x-y\|}\le b_{n_k}(\gamma_2\gamma_3)^2$ for all $x\ne y\in X$ and $(I,J)\in{\mathcal B}_k$, where $\gamma_2>0$ is defined in \eqref{eq3.5} and $\gamma_3$ is defined in  \eqref{eq3.6} when $\rho$ is replaced by the Euclidean norm.
Therefore, for the $X_1$ in Definition \ref{def2.2} (ib), we have $\phi_{1,I}(X)\subset\phi_{1,J}(X_1)$ for all large $k>0$ and $(I,J)\in{\mathcal B}_k$. Hence, the assumptions in Definition \ref{def2.2} (ib) show that $\phi_{1,J}^{-1}$ is well defined in $\phi_{1,I}(X)$ for all large $k>0$ and $(I,J)\in{\mathcal B}_k$ and,
\[
\lim_{k\to\infty}\sup_{(I, J)\in{\mathcal B}_k}\sup_{x\in X}\|\phi_{1,J}^{-1}\circ\phi_{1,I}(x)-x\|=0,\quad \lim_{k\to\infty}\#{\mathcal B}_k=+\infty.
\]

If $\{\Phi_n\}_n$ is a self-affine MIFS, by the MWHP and the definition of $S_J$, we see that
$\gamma_2^{-2}\le\frac{\|S_{J}(x)-S_{J}(y)\|}{\|x-y\|}\le \gamma_2^2$ and $\gamma_2^{-2}\le\frac{\|S_{I}(x)-S_{I}(y)\|}{\|x-y\|}\le \gamma_2^2$ for all $x\ne y\in X$, where $\gamma_2>0$ is defined in \eqref{eq3.5}.
Recall that all $\phi_{n,j}$ can be extended to ${\mathbb R}^d$. It is easy to see that
\[
\lim_{k\to\infty}\sup_{(I, J)\in{\mathcal B}_k}\sup_{x\in X}\|\phi_{1,J}^{-1}\circ\phi_{1,I}(x)-x\|=0,\quad \lim_{k\to\infty}\#{\mathcal B}_k=+\infty.
\]
The necessity of (i) is proven.

\medskip

Finally, we prove the sufficiency of (i).

\medskip

Suppose there exists a sequence of numbers $b_n\in (0,\ 1)$, a sequence of words $\sigma_n\in{\mathcal I}_{b_n}$ and a sequence of pairs of words ${\mathcal B}_n$ such that \eqref{eqww8-23} holds. It is   an easily observable fact that either the number of all different $I$ satisfying $(I,J)\in{\mathcal B}_n$ for some $J\in{\mathcal I}_{\sigma_n}$ is at least $\sqrt{\#{\mathcal B}_n}$ or the number of all different $J$ satisfying $(I,J)\in{\mathcal B}_n$ for some $I\in{\mathcal I}_{\sigma_n}$ is at least $\sqrt{\#{\mathcal B}_k}$. Without loss of generalization, assume that the number of all different $J$ satisfying $(I,J)\in{\mathcal B}_n$ for some $I\in{\mathcal I}_{\sigma_n}$ is at least $\sqrt{\#{\mathcal B}_k}$. Hence, the inclusion relationship in \eqref{eqww8-23} shows $$\#\left\{J\in{\mathcal I}_{\sigma_n}:\ \phi_{1,J}(X)\cap\phi_{1,\sigma_n}(X)\ne\emptyset\right\}\ge \sqrt{\#{\mathcal B}_n}.$$ Therefore, we get
\[\begin{array}{rl}
&\sup\limits_{1>b>0}\max\limits_{I\in{\mathcal I}_b}\#\left\{J\in{\mathcal I}_b:\ \phi_{1,J}(X)\cap\phi_{1,I}(X)\ne\emptyset\right\}\medskip\\
\ge& \lim\limits_{n\to\infty}\#\left\{J\in{\mathcal I}_{b_n}:\ \phi_{1,J}(X)\cap\phi_{1,\sigma_n}(X)\ne\emptyset\right\}\medskip\\
\ge &\lim\limits_{n\to\infty}\sqrt{\#{\mathcal B}_n}=+\infty.
\end{array}\]
According to Lemma \ref{th3.5}, we conclude that the IMFS can not satisfy  MOSC. This proves the sufficiency of (i).

\medskip

(ii) The proof is similar to (i) if we replace $J$ by $\phi_{1,J}$.
\end{pf}


\section{Dimensions of Moran-type self-similar sets}
\setcounter{equation}{0}\setcounter{theorem}{0}
In this section, we first consider the box dimension and the packing dimension.

\begin{theorem}\label{th4.1} Let $\left\{\Phi_n\right\}_{n=1}^\infty$ be a self-similar MIFS defined in Definition \ref{def2.2} (i) where $\phi_{n,j}(x)=r_{n,j}O_{n,j}(x+\alpha_{n,j})$, $r_{n,j}\in(0,\ 1)$ and $O_{n,j}$ are orthogonal matrices. Assume the MIFS satisfies the MWSC defined in Definition \ref{def3.1}.

\medskip

{\rm (i)} If $r_0:=\inf\{r_{n,j}:\ 1\le j\le N_n,\ n\ge 1\}>0$, then
$\underline{\dim}_{\rm B}K_n=\liminf\limits_{b\to0}\frac{\ln\#{\mathcal A}_b}{-\ln b}\leq \limsup\limits_{b\to0}\frac{\ln\#{\mathcal A}_b}{-\ln b}=\overline{\dim}_{\rm B}K_n=\dim_{\rm P}K_n$
for all $n\ge1$.

\medskip

{\rm (ii)} If $r_{n,j}=r_n$ for any  $n>0$ and $1\le j\le N_n$, and

\begin{equation}\label{contra-ra}
\lim\limits_{n\to+\infty}\frac{\ln r_n}{\ln (r_1r_2\cdots r_{n-1})}=0,
\end{equation}
 then, for any $n\geq1$ we have
\begin{equation}\label{lower-dim}
\underline{\dim}_{\rm B}K_n=\liminf\limits_{b\to0}\frac{\ln\#{\mathcal A}_b}{-\ln b}=\liminf\limits_{n\to+\infty}\frac{\ln\#{\mathcal A}_{r_1r_2\cdots r_n}}{-\ln (r_1r_2\cdots r_n)},
\end{equation}
and \begin{equation}\label{upper-dim}
\overline{\dim}_{\rm B}K_n=\dim_{\rm P}K_n= \limsup\limits_{b\to0}\frac{\ln\#{\mathcal A}_b}{-\ln b}=\limsup\limits_{n\to+\infty}\frac{\ln\#{\mathcal A}_{r_1r_2\cdots r_n}}{-\ln (r_1r_2\cdots r_n)}.
\end{equation}
\end{theorem}
\begin{pf} By equalities in \eqref{eq2.2} of Theorem \ref{eq2.1}, we see that $\underline{\dim}_{\rm B}K_n$, $\overline{\dim}_{\rm B}K_n$ and $\dim_{\rm P}K_n$ are independent of $n>0$. Thus, we need only to consider the dimensions of $K_1$.

\medskip

For any open set $U\subset{\mathbb R}^d$, if $U\cap K_1\ne\emptyset$, then \eqref{eq1.1} shows that there exists an integer $n>0$ and a word $J\in\Sigma_1^n$ such that $\phi_{1,J}(K_{n+2})\subset U$. Hence, \cite[Corollary 3.9]{[falconer]} shows $\overline{\dim}_{\rm B}K_1=\dim_{\rm P}K_1$ by using $\overline{\dim}_{\rm B}K_1=\overline{\dim}_{\rm B}K_n$ and $\dim_{\rm P}K_1=\dim_{\rm P}K_n$. Therefore, we need only to consider $\underline{\dim}_{\rm B}K_1$ and $\overline{\dim}_{\rm B}K_1$ in the following.

\medskip

(i) For any integer vector ${\bf z}=(z_1,\cdots,z_d)^t\in{\mathbb Z}^d$ and positive number $b>0$, denote the cube $\prod_{j=1}^d[z_jb|X|,\ (z_j+1)b|X|]$ by ${\mathcal C}_{b,{\bf z}}$. Let
\[
N_b=\#\{{\mathcal C}_{b,{\bf z}}:\ {\bf z}\in{\mathbb Z}^d,\ K_1\cap {\mathcal C}_{b,{\bf z}}\ne\emptyset\}.
\]
It is well known that $\underline{\dim}_{\rm B}K_1=\liminf\limits_{b\to0^+}\frac{\ln N_b}{-\ln b}$ and $\overline{\dim}_{\rm B}K_1=\limsup\limits_{b\to0}\frac{\ln N_b}{-\ln b}$.

\medskip

 From the relationship $K_1\subseteq\bigcup\limits_{f\in{\mathcal A}_b}f(X)$, it follows that if ${\mathcal C}_{b,{\bf z}}\cap K_1\ne \emptyset$, then  ${\mathcal C}_{b,{\bf z}}$ intersects at least one set $f(X)$ with $f\in{\mathcal A}_b$. Consequently, we have
\[
\{{\mathcal C}_{b,{\bf z}}:\ {\bf z}\in{\mathbb Z}^d,\ K_1\cap {\mathcal C}_{b,{\bf z}}\ne\emptyset\}\subset\bigcup_{f\in{\mathcal A}_b}\{{\mathcal C}_{b,{\bf z}}:\ {\bf z}\in{\mathbb Z}^d,\ f(X)\cap {\mathcal C}_{b,{\bf z}}\ne\emptyset\}.
\]
Secondly, since $|f(X)|\le b|X|$ for any $f\in {\mathcal A}_b$, we see that $f(X)$ intersects no more than $3^d$ cubes for any $f\in {\mathcal A}_b$. Hence, $\#\{{\mathcal C}_{b,{\bf z}}:\ {\bf z}\in{\mathbb Z}^d,\ f(X)\cap {\mathcal C}_{b,{\bf z}}\ne\emptyset\}\leq 3^d$ for any $f\in {\mathcal A}_b$. As a result, we have
\begin{equation}\label{left-ineq}
N_b\leq 3^d\#{\mathcal A}_b.
\end{equation}

\medskip

As  $r_0=\inf\{r_{n,j}:\ 1\le j\le N_n,\ n\ge 1\}>0$, the MIFS $\{\Phi_n\}_{n\geq1}$ satisfies the MWHP. Fix a cube ${\mathcal C}_{b,{\bf z}}$, from Lemma  \ref{th3.5}, we see  that there is a sub family $\Psi_{b,{\bf z}}\subseteq\{f\in{\mathcal A}_b:\ f(X)\cap{\mathcal C}_{b,{\bf z}}\ne\emptyset\}$ such that: (a) $\gamma_4\#\Psi_{b,{\bf z}}\ge \#\{f\in{\mathcal A}_b:\ f(X)\cap{\mathcal C}_{b,{\bf z}}{\neq\emptyset}\}$, where $\gamma_4$ is defined in \eqref{eq3.7}; (b) Every $f\in \{f\in{\mathcal A}_b:\ f(X)\cap{\mathcal C}_{b,{\bf z}}\ne\emptyset\}$ corresponds at least one $g\in\Psi_{b,{\bf z}}$ such that $f(X)\cap g(X)\ne\emptyset$. (c) $f(X)\cap h(X)=\emptyset$ for any $f\ne h\in\Psi_{b,{\bf z}}$.
Then, the definition of ${\mathcal A}_b$ shows that $\bigcup_{f\in\Psi_{b,{\bf z}}}f(X)$ is a disjoint union contained in a union of $3^d$ cubes since $|f(X)|\le b|X|$. For any $ f=\phi_{1,j_1j_2\cdots j_n}\in{\mathcal A}_b$,  by computing Lebesgue's measure, we have
\[
{\mathcal L}^d(f(X))=(r_{1,j_1j_2\cdots j_{n-1}}r_{n,j_n})^d{\mathcal L}^d(X)\ge (br_{n,j_n})^d{\mathcal L}^d(X)\ge (br_0)^d{\mathcal L}^d(X).
\]
 Hence, we have
\[
\begin{array}{rl}
\#\{f\in{\mathcal A}_b:\ f(X)\cap{\mathcal C}_{b,{\bf z}}\ne\emptyset\}\leq &\gamma_4\#\Psi_{b,{\bf z}}\\ \leq & \gamma_4(br_0)^{-d}({\mathcal L}^d(X))^{-1}{\mathcal L}^d\left(\bigcup\limits_{f\in\Psi_{b,{\bf z}}}f(X)\right)\\ \leq &\gamma_4 r_0^{-d}({\mathcal L}^d(X))^{-1}(3|X|)^d.
\end{array}
\]
Therefore,
\[
 \ \ \#{\mathcal A}_b\leq \gamma_4 r_0^{-d}({\mathcal L}^d(X))^{-1}(3|X|)^dN_b,\quad \forall\ b\in(0,\ 1).
\]
Together with \eqref{left-ineq}, it leads to  that there is a constant $c\geq 1$ such that
\begin{equation}\label{eq4.1}
c^{-1}N_b\leq \#{\mathcal A}_b\leq cN_b,\quad \forall\ b\in(0,\ 1).
\end{equation}
Hence, we have
 \[\underline{\dim}_{\rm B}K_1=\liminf\limits_{b\to0^+}\frac{\ln N_b}{-\ln b}=\liminf\limits_{b\to0^+}\frac{\ln\#{\mathcal A}_b}{-\ln b}\leq\limsup\limits_{b\to0^+}\frac{\ln\#{\mathcal A}_b}{-\ln b}=\limsup\limits_{b\to0^+}\frac{\ln N_b}{-\ln b}=\overline{\dim}_{\rm B}K_1.\]

\medskip

(ii) If $\inf\{r_n:\ n>0\}>0$, the conclusion follows from (i). We  turn to the case $\inf\{r_n:\ n>0\}=0$ in the following.

\medskip

It is easy to see that \eqref{left-ineq} still holds since the inequality is independent of the assumption $r_0>0$.

\medskip

As  $r_{n,j}=r_n$ for any  $n>0$ and $1\le j\le N_n$, one see that for any $b\in (0,r_1)$, there exists $n_b>0$ such that ${\mathcal I}_b=\Sigma_1^{n_b}$ and $r_1r_2\cdots r_{n_b}\leq b<r_1r_2\cdots r_{n_b-1}$. Hence, we have
\[
\frac{\ln\#{\mathcal A}_{r_1r_2\cdots r_{n_b}}}{-\ln (r_1r_2\cdots r_{n_b-1})}\geq\frac{\ln\#{\mathcal A}_b}{-\ln b}\geq\frac{\ln\#{\mathcal A}_{r_1r_2\cdots r_{n_b}}}{-\ln (r_1r_2\cdots r_{n_b})}.
\]
In combination with \eqref{contra-ra}, we have
\begin{equation}\label{eq4.2}
\left\{\begin{array}{l}\liminf\limits_{b\to0^+}\frac{\ln\#{\mathcal A}_b}{-\ln b}=\liminf\limits_{n\to+\infty}\frac{\ln\#{\mathcal A}_{r_1r_2\cdots r_n}}{-\ln (r_1r_2\cdots r_n)},\\ \limsup\limits_{b\to0^+}\frac{\ln\#{\mathcal A}_b}{-\ln b}=\limsup\limits_{n\to+\infty}\frac{\ln\#{\mathcal A}_{r_1r_2\cdots r_n}}{-\ln (r_1r_2\cdots r_n).}
\end{array}\right.\end{equation}
As all assumptions in Lemma \ref{th3.5} are satisfied, almost the same argument of \eqref{eq4.1} shows that there is a constant $c>0$ such that
\begin{equation}\label{eq4.3}
c^{-1}N_{r_1r_2\cdots r_n}\leq \#{\mathcal A}_{r_1r_2\cdots r_n}\leq cN_{r_1r_2\cdots r_n},\quad  \forall n\geq1.
\end{equation}
In combination with \eqref{left-ineq} and  \eqref{eq4.2}, the above inequality implies
\begin{equation}\label{eq4.4}
\begin{array}{rl}
\liminf\limits_{b\to0^+}\frac{\ln N_b}{-\ln b}\leq\liminf\limits_{b\to0^+}\frac{\ln\#{\mathcal A}_{b}}{-\ln b}
\le\liminf\limits_{n\to+\infty}\frac{\ln\#{\mathcal A}_{r_1r_2\cdots r_n}}{-\ln (r_1r_2\cdots r_n)}=\liminf\limits_{n\to+\infty}\frac{\ln N_{r_1r_2\cdots r_n}}{-\ln (r_1r_2\cdots r_n)}.
\end{array}
\end{equation}
By a well-known technique (p.$44$,\ \cite{ [falconer]}), we have
 \[\liminf\limits_{b\to0^+}\frac{\ln N_b}{-\ln b}=\liminf\limits_{n\to+\infty}\frac{\ln N_{2^{-n}}}{n\ln 2}.\]
  For any integer $k>-\log_2r_1$, there is a unique integer $n_k>0$ such that
\[
r_1r_2\cdots r_{n_k}\leq 2^{-k}<r_1r_2\cdots r_{n_k-1},\quad k>0.
\]
Hence, there exist integers $s_k$ and $\ t_k$ with $s_k\leq k\leq t_k$ such that
\[\left\{\begin{array}{l}
2^{-t_k-1}< r_1r_2\cdots r_{n_k}\leq 2^{-t_k}\leq 2^{-k},\\ 2^{-k}\leq 2^{-s_k}\leq r_1r_2\cdots r_{n_k-1}<2^{-s_k+1},
\end{array}\right.
\]
for all $ k>-\log_2r_1$. Noting that $N_{2^{-n}}$ increases as $n$ increases, we have
\[
\frac{\ln N_{2^{-k}}}{k\ln 2}\ge \frac{\ln N_{2^{-s_k}}}{t_k\ln 2}.
\]
Due to \eqref{contra-ra}, it follows that $\lim\limits_{k\to+\infty}\frac{s_k}{t_k}=1$. Therefore, the property of lower limit of sequences shows

\begin{equation}
\begin{array}{rl}
\liminf\limits_{b\to0^+}\frac{\ln N_b}{-\ln b}=&\liminf\limits_{n\to+\infty}\frac{\ln N_{2^{-n}}}{n\ln 2}\geq \liminf\limits_{k\to+\infty}\frac{\ln N_{2^{-s_k}}}{t_k\ln 2}\\= &\liminf\limits_{k\to+\infty}\frac{\ln N_{2^{-s_k}}}{s_k\ln 2}\geq \liminf\limits_{k\to+\infty}\frac{\ln N_{r_1r_2\cdots r_{n_k-1}}}{-\ln (r_1r_2\cdots r_{n_k-1})}\\
\geq &\liminf\limits_{n\to+\infty}\frac{\ln N_{r_1r_2\cdots r_n}}{-\ln (r_1r_2\cdots r_n)}\geq \liminf\limits_{b\to0^+}\frac{\ln N_b}{-\ln b}.
\end{array}
\end{equation}
Together with \eqref{eq4.4}, the above inequality leads to
 \[\underline{\dim}_B(K_1)=\liminf\limits_{b\to0^+}\frac{\ln N_b}{-\ln b}=\liminf\limits_{b\to0^+}\frac{\ln\#{\mathcal A}_b}{-\ln b}=\liminf\limits_{n\to+\infty}\frac{\ln \#{\mathcal A}_{r_1r_2\cdots r_n}}{-\ln (r_1r_2\cdots r_n)}.\]

Next, we turn to $\overline{\dim}_B(K_1)$.
 By \eqref{left-ineq}, we have

\begin{equation}\label{sup-ineq}
\limsup\limits_{b\to0^+}\frac{\ln N_b}{-\ln b}\leq\limsup\limits_{b\to0^+}\frac{\ln\#{\mathcal A}_{b}}{-\ln b}.
\end{equation}
In combination \eqref{eq4.2} and \eqref{eq4.3}, we have
\[
\limsup\limits_{b\to0^+}\frac{\ln N_b}{-\ln b}\geq \limsup\limits_{n\to+\infty}\frac{\ln N_{r_1r_2\cdots r_n}}{-\ln (r_1r_2\cdots r_n)}=\limsup\limits_{n\to+\infty}\frac{\ln\#{\mathcal A}_{r_1r_2\cdots r_n}}{-\ln (r_1r_2\cdots r_n)}=\limsup\limits_{b\to0^+}\frac{\ln\#{\mathcal A}_b}{-\ln b}.
\]
Together with \eqref{sup-ineq}, the above inequalities show
\[\overline{\dim}_B(K_1)=\limsup\limits_{b\to0^+}\frac{\ln N_b}{-\ln b}=\limsup\limits_{b\to0^+}\frac{\ln\#{\mathcal A}_b}{-\ln b}=\limsup\limits_{n\to+\infty}\frac{\ln \#{\mathcal A}_{r_1r_2\cdots r_n}}{-\ln (r_1r_2\cdots r_n)}.\]

\medskip

The proof is finished.
\end{pf}

Next, we turn to  the Hausdorff dimension of $K_n$.
\medskip

\begin{theorem}\label{th4.3} Let $\left\{\Phi_n\right\}_{n=1}^\infty$ be a self-similar MIFS defined in Definition \ref{def2.2} (i) where $\phi_{n,j}(x)=r_{n,j}O_{n,j}(x+\alpha_{n,j})$, $r_{n,j}\in(0,\ 1)$ and $O_{n,j}$ are orthogonal matrices. Assume that the MIFS satisfies the MOSC defined in Definition \ref{def3.1} and $r_0:=\inf\{r_{n,j}:\ 1\le j\le N_n,\ n\ge 1\}>0$.

\medskip

\item[\textup{(i)}] If  $+\infty>\liminf\limits_{n\to\infty}\prod\limits_{i=1}^n\sum\limits_{j=1}^{N_i}r_{i,j}^s>0$, then $+\infty>{\mathcal H}^s(K_1)>0$ and so $\dim_H(K_1)=s$.

\medskip

\item[\textup{(ii)}] If $\liminf\limits_{n\to\infty}\prod\limits_{i=1}^n\sum\limits_{j=1}^{N_i}r_{i,j}^s=+\infty$, then ${\mathcal H}^s(K_1)=+\infty$ and so $\dim_H(K_1)\ge s$.

\medskip

\item[\textup{(iii)}] If $\liminf\limits_{n\to\infty}\prod\limits_{i=1}^n\sum\limits_{j=1}^{N_i}r_{i,j}^s=0$, then ${\mathcal H}^s(K_1)=0$ and so $\dim_H(K_1)\le s$.
\end{theorem}
\begin{pf}
(i)
Write $\alpha=\liminf\limits_{n\to\infty}\prod\limits_{i=1}^n\sum\limits_{j=1}^{N_i}r_{i,j}^s.$ There exists a positive integer $N$ such that
\begin{equation}\label{product}
\prod\limits_{i=1}^n\sum\limits_{j=1}^{N_i}r_{i,j}^s>\frac{\alpha}{2}>0,\quad\forall\ n>N.
\end{equation}
Let
\[
p_{n,j}=\frac{\gamma_{n,j}^s}{\sum_{j=1}^{N_i}\gamma_{i,j}^s},\quad\forall\ 1\le j\le N_n,\ n\ge1.
\]
Then, there is a probability measure on $\Sigma_1^\mathbb N$ such that
\[
\nu(\Sigma_1^{\mathbb N}(j_1,\cdots, j_n))=p_{1,j_1}\cdots p_{n,j_n}.
\]

In combination with \eqref{product}, for any $n\geq1$ large enough, we have

\begin{equation}\label{eq4.11}
\begin{array}{rl}
\nu(\Sigma_1^{\mathbb N}(j_1j_2\cdots j_{n}))=&p_{1,j_1}p_{2,j_2}\cdots p_{n,j_n}\medskip\\
=& \left(\prod\limits_{i=1}^{n}[\sum\limits_{j=1}^{N_i}r_{i,j}^s]\right)^{-1}(r_{1,i_1}r_{2,i_2}\cdots r_{n,i_n})^s\\
\leq& \frac{2}{\alpha}(r_{1,j_1}r_{2,j_2}\cdots r_{n,j_n})^s.
\end{array}
\end{equation}

Note that \eqref{eq2.1} defines a map $\pi:\ \Sigma_1^{\mathbb N}\to K_1$ for any fixed $a\in X$. By \eqref{eq1.1} and the compactness of $X$, we see that this $\pi$ is a continuous map from $(\Sigma_1^{\mathbb N},\ \varrho_1)$ to $(K_1,\ \rho)$. Hence $\nu\circ\pi^{-1}$ is a Borel probability measure on $K_1$.

\medskip

In virtue of  $r_0>0$, we see that the self-similar MIFS $\left\{\Phi_n\right\}_{n=1}^\infty$ satisfies  the MWHP and MBDP. By a  similar argument  to Lemma \ref{th3.5}, we have
\begin{equation}\label{eq4.13}
\gamma_4':=\sup\limits_{U\subset X}\#\left\{{J}\in{\mathcal I}_{|U|}:\ \phi_{1,J}(X)\cap U\ne\emptyset\right\}<+\infty.
\end{equation}
As a consequence, we get $\#\{i_1i_2\cdots i_{n}\in{\mathcal I}_{|U|}:\ \phi_{1,i_1i_2\cdots i_{n}}(K_{n+1})\cap U\ne\emptyset\}\le\gamma_4'$. This means that $\pi^{-1}(U)$ is contained in a union of at most $\gamma'_4$ cylinder sets of the form $\Sigma_1^{\mathbb N}(i_1i_2\cdots i_{n})$ with $i_1i_2\cdots i_{n}\in{\mathcal I}_{|U|}$. Hence, there is a finite word $i_1i_2\cdots i_{n}\in{\mathcal I}_{|U|}$ such that
$\nu\circ\pi^{-1}(U)\leq \gamma'_4\nu(\Sigma_1^{\mathbb{N}}(i_1\cdots i_n))$. Consequently,  we get

\begin{equation}\label{eq4.14}
\nu\circ\pi^{-1}(U)\leq \frac{2}{\alpha}\gamma_4' (r_{1,i_1}r_{2,i_2}\cdots r_{n,i_n})^s
\leq \frac{2}{\alpha}\gamma_4' |U|^s
\end{equation}
by using \eqref{eq4.11}. Therefore, the famous mass distribution principle shows ${\mathcal H}^s(K_1) \ge \frac{\alpha}{2\gamma_4'}>0$, and so $\dim_H(K_1)\geq s$.

\medskip

Finally, we prove ${\mathcal H}^s(K_1)<+\infty$. Since $\bigcup\limits_{J\in\Sigma_1^n}\phi_{1,J}(X)$ is a $\delta_n$-cover of $K_1$ with $\delta_n=\max\{|\phi_{1,J}(X)|:\ J\in\Sigma_1^n\}$, the definition of Hausdorff measures shows
\[
{\mathcal H}^s(K_1)\le\liminf\limits_{n\to\infty}\sum\limits_{J\in\Sigma_1^n}|\phi_{1,J}(X)|^s
=\liminf\limits_{k\to\infty}\sum\limits_{J\in\Sigma_1^n}r_{1,J}^s|X|^s<+\infty
\]
by noting
$\liminf\limits_{n\to\infty}\sum\limits_{J\in\Sigma_1^{n}}r_{1,J}^s=\liminf\limits_{n\to\infty}\prod\limits_{i=1}^n\sum\limits_{j=1}^{N_i}r_{i,j}^s<+\infty$. This implies  that ${\mathcal H}^s(K_1)<+\infty$
and $\dim_{H}(K_1)\le s$.
Therefore, conclusion (i) is proven.

\medskip

(ii) Assume that $\liminf\limits_{n\to\infty}\prod\limits_{i=1}^n\sum\limits_{j=1}^{N_i}r_{i,j}^s=+\infty$. By \eqref{eq4.11}, we can rewrite  \eqref{eq4.14}  as
\[
\nu\circ\pi^{-1}(U)\leq c_{|U|}\gamma_4' |U|^s
\]
with $\lim\limits_{|U|\to0}c_{|U|}=0$. Then, the definition of Hausdorff measures implies
\[\begin{array}{rl}
{\mathcal H}^s(K_1)=&\lim\limits_{\delta\to0}\inf\{\sum\limits_{j}|U_j|^s:\ K_1 \subseteq \bigcup U_j,\ |U_j|\leq \delta\}\\
\ge &\lim\limits_{\delta\to0}\inf\{\sum\limits_{j}(c_{|U_j|}\gamma_4' {}p)^{-1}\nu\circ\pi^{-1}(U_j):\ K_1 \subseteq \bigcup U_j,\ |U_j|\leq \delta\}\\
\ge &\lim\limits_{\delta\to0}\inf\{\min\limits_{j}(c_{|U_j|}\gamma_4' {})^{-1}:\ K_1 \subseteq \bigcup U_j,\ |U_j|\leq \delta\}\\
=&+\infty,
\end{array}\]
where the last equality follows from $|U_j|\leq \delta$ and $\lim\limits_{|U|\to0}c_{|U|}=0$.
   Hence, the conclusion (ii) is proven.

\medskip

(iii) The proof is almost the same as the second part of the proof of (i).
\end{pf}

\medskip

\begin{coro}\label{th4.4} Let $\left\{\Phi_n\right\}_{n=1}^\infty$ be a self-similar MIFS defined in Definition \ref{def2.2} (i), where $\phi_{n,j}(x)=r_{n,j}O_{n,j}(x+\alpha_{n,j})$, $r_{n,j}\in(0,\ 1)$ and $O_{n,j}$ are orthogonal matrices. Assume the MIFS satisfies the MOSC defined in Definition \ref{def3.1} and $r_0:=\inf\{r_{n,j}:\ 1\le j\le N_n,\ n\ge 1\}>0$. Let $s_k$ be the unique solution of $\prod\limits_{i=1}^k\sum\limits_{j=1}^{N_i}r_{i,j}^s=1$ for each $k\ge1$. The following statements hold.

\medskip

{\rm (i)} $\dim_H(K_1)=\liminf\limits_{k\to+\infty} s_k$.

\medskip

{\rm (ii)} Let $\dim_H(K_1)=s$. Then, we have
\[
\left\{\begin{array}{ll}
{\mathcal H}^s(K_1)=+\infty,&\ \ {\rm if}\ \liminf\limits_{n\to\infty}\prod\limits_{i=1}^n\sum\limits_{j=1}^{N_i}r_{i,j}^s=+\infty,\\
{\mathcal H}^s(K_1)=0,&\ \ {\rm if}\ \liminf\limits_{n\to\infty}\prod\limits_{i=1}^n\sum\limits_{j=1}^{N_i}r_{i,j}^s=0,\\
0<{\mathcal H}^s(K_1)<+\infty,&\ \ {\rm if}\ 0<\liminf\limits_{n\to\infty}\prod\limits_{i=1}^n\sum\limits_{j=1}^{N_i}r_{i,j}^s<+\infty.
\end{array}\right.
\]

\medskip

{\rm (iii)} If $r_{n,j}=r_n$ for all $n>0$ and $1\le j\le N_n$, then $\dim_H(K_1)=\underline{\dim}_B(K_1)$.
\end{coro}
\begin{pf}
\medskip
(i)
Assume that  $t<\liminf\limits_{k\to+\infty} s_k$. There exist an integer $k_0>0$ and a positive number $\varepsilon>0$ such that $s_k>t+\varepsilon$ for all $k\ge k_0$. In virtue of \eqref{eq1.1} and the definition of $s_k$, we have

\[\prod\limits_{i=1}^n\sum\limits_{j=1}^{N_i}r_{i,j}^t>\prod\limits_{i=1}^n\sum\limits_{j=1}^{N_i}r_{i,j}^{s_n-\varepsilon}\ge \prod\limits_{i=1}^n\sum\limits_{j=1}^{N_i}r_{i,j}^{s_n}c_{2,j}^{-\varepsilon}=\prod\limits_{i=1}^n\sum\limits_{j=1}^{N_i}c_{2,j}^{-\varepsilon} \to+\infty,\ n\to +\infty,\]
which implies that $\liminf\limits_{n\to\infty}\prod\limits_{i=1}^n\sum\limits_{j=1}^{N_i}r_{i,j}^t=+\infty$. As a consequence of  Theorem \ref{th4.3} (ii), we get ${\mathcal H}^t(K_1)=+\infty$, and so $\dim_H(K_1)\ge t$. By the arbitrariness of  $t$, we have $\dim_H(K_1)\ge \liminf\limits_{k\to+\infty} s_k$.

\medskip

Let $t>\liminf\limits_{k\to+\infty} s_k$. There exist a sequence of integers $n_\ell$ and a positive number $\varepsilon>0$ such that $s_{n_\ell}<t-\varepsilon$ for all $\ell>0$ and $\lim\limits_{\ell\to+\infty} s_{n_\ell}=\liminf\limits_{k\to+\infty} s_k$.
From \eqref{eq1.1}, it follows that
\[\prod\limits_{i=1}^{n_\ell}\sum\limits_{j=1}^{N_i}r_{i,j}^t\leq \prod\limits_{i=1}^{n_\ell}\sum\limits_{j=1}^{N_i}r_{i,j}^{s_{n_\ell}+\varepsilon}\le \prod\limits_{i=1}^{n_\ell}\sum\limits_{j=1}^{N_i}r_{i,j}^{s_{n_\ell}}c_{2,j}^\varepsilon= \prod\limits_{i=1}^{n_\ell}c_{2,j}^\varepsilon
\to 0,\ \ell\to +\infty.
\]
 It follows that $\liminf\limits_{n\to+\infty}\prod\limits_{i=1}^n\sum\limits_{j=1}^{N_i}r_{i,j}^{t}=0$. By Theorem \ref{th4.3} (iii), we see ${\mathcal H}^t(K_1)=0$, and so $\dim_H(K_1)\le t$. Therefore, $\dim_H(K_1)\le \liminf\limits_{k\to+\infty} s_k$.

\medskip

In conclusion, $\dim_H(K_1)=\liminf\limits_{k\to+\infty} s_k$.

\medskip

(ii) It  directly follows from Theorem \ref{th4.3}.

\medskip

(iii) It follows from combining  Theorem \ref{th4.1} (i), the conclusion (i) and the fact that the assumption $r_{n,j}=r_n$ implies ${\mathcal A}_b=\{\phi_{1,J}:\ J\in\Sigma_1^n\}$ when $r_1r_2\cdots r_n\leq b<r_1r_2\cdots r_{n-1}$.
\end{pf}

\medskip

For the MIFS, it is interesting whether $K_n$ is an $s$-set. The following result include the answer and the relationship between Hausdorff dimension and box dimension.

\begin{theorem}\label{th4.5} Let $\left\{\Phi_n\right\}_{n=1}^\infty$ be a self-similar MIFS defined in Definition \ref{def2.2} (i), where $\phi_{n,j}(x)=r_{n,j}O_{n,j}(x+\alpha_{n,j})$, $r_{n,j}\in(0,\ 1)$, $O_{n,j}$ are orthogonal matrices. Assume the MIFS satisfies the MOSC defined in Definition \ref{def3.1}.
If $0<\liminf\limits_{n\to\infty}\prod\limits_{i=1}^n\sum\limits_{j=1}^{N_i}r_{i,j}^s$ $\le\limsup\limits_{n\to\infty}\prod\limits_{i=1}^n\sum\limits_{j=1}^{N_i}r_{i,j}^s<+\infty$ and $r_0:=\inf\{r_{n,j}:\ n>0,\ 1\leq j\le N_n\}>0$, then the following statements hold.

\medskip

{\rm (i)} All $K_{k}$ are uniform $s$-sets:
\begin{equation}\label{eq4.15}
+\infty>\sup\{{\mathcal H}^s(K_{k}):\ k>0\}\ge\inf\{{\mathcal H}^s(K_{k}):\ k>0\}>0.
\end{equation}

{\rm (ii)} $\dim_{H}K_n=\underline{\dim}_{\rm B}K_n=\overline{\dim}_{\rm B}K_n=s$ for all $n>0$.
\end{theorem}
\begin{pf}
(i) The assumption $0<\liminf\limits_{n\to\infty}\prod\limits_{i=1}^n\sum\limits_{j=1}^{N_i}r_{i,j}^s\le\limsup\limits_{n\to\infty}\prod\limits_{i=1}^n\sum\limits_{j=1}^{N_i}r_{i,j}^s<+\infty$ implies that there exist two positive numbers $a\le b$ such that
\[
a\leq\prod\limits_{i=1}^n\sum\limits_{j=1}^{N_i}r_{i,j}^s\le b,\quad \forall\ n>0,
\]
so
\[
ab^{-1}\leq\prod\limits_{i=m+1}^n\sum\limits_{j=1}^{N_i}r_{i,j}^s\le ba^{-1},\quad \forall\ n>m\ge1.
\]
Therefore, we get
\[
a^2b^{-1}\leq\prod\limits_{i=m}^n\sum\limits_{j=1}^{N_i}r_{i,j}^s\le b^2a^{-1},\quad \forall\ n>m\ge1,
\]
so
\[
a^2b^{-1}\leq\liminf\limits_{n\to\infty}\prod\limits_{i=m}^n\sum\limits_{j=1}^{N_i}r_{i,j}^s\le\limsup\limits_{n\to\infty}\prod\limits_{i=m}^n\sum\limits_{j=1}^{N_i}r_{i,j}^s\le b^2a^{-1},\quad \forall\ m\ge1.
\]
Hence, the definition of Hausdorff measures implies
\begin{equation}\label{eq4.16}
{\mathcal H}^s(K_{m})\le\liminf\limits_{k\to\infty}\sum\limits_{J\in\Sigma_{m}^{k}}|\phi_{m,J}(X)|^s\leq b^2a^{-1}|X|^s,\quad \forall\ m>0.
\end{equation}

On the other hand, consider the sub MIFS $\left\{\Phi_n\right\}_{n=m+1}^\infty$ for any given $m>0$. 
By a similar argument to \eqref{eq4.11}, we have
\[
\nu(\Sigma_{m+1}^{\mathbb N}(i_{m+1}i_{m+2}\cdots i_{n}))\leq ba^{-1}(r_{m+1,i_{m+1}}r_{m+2,i_{m+2}}\cdots r_{n,i_n})^s.
\]
As a result, we also get
\[
\nu\circ\pi^{-1}(U)\leq  ba^{-1}\gamma_4'{ }|U|^s.
\]
Note that the definition of $\gamma_4'$ shows that it is non-increasing when $m$ increases.  Then, it can be chosen to be independent of $m$. Therefore,
\begin{equation}\label{lowerbound}
\inf\{{\mathcal H}^s(K_{m}):\ m\ge1\}\ge(ba^{-1}\gamma_4' {})^{-1}>0
\end{equation}
by using the mass distribution principle. Hence, \eqref{eq4.15} holds by using this inequality and \eqref{eq4.16}. We finish the proof of (i)

\medskip

(ii) Since our MIFS satisfies the MOSC and the assumption $r_0>0$ implies that our MIFS satisfies the MWHP, from Lemma \ref{th3.5} (ii) it follows that for every $b\in(0,\ 1)$, there exists a sub family $\{\phi_{1,J_1},\ \phi_{1,J_2},\ \dots,\ \phi_{1,J_{k_b}}\}\subset{\mathcal A}_b$ such that
(a) $\phi_{1,J_1}(X),\dots,\phi_{1,J_{k_b}}(X)$ are pairwise disjoint, (b) $\inf\{\frac{k_b}{\#{\mathcal A}_b}:\ 0<b<1\}>0$.
By \eqref{lowerbound}, we have
\[
\begin{array}{rl}
{\mathcal H}^s(K_{1})\geq\sum\limits_{i=1}^{k_b}r_{1,J_i}^s{\mathcal H}^s(K_{|J_i|+1}) \geq k_b b^sr^s_0  \frac{a}{b\gamma_4'} \geq \#{\mathcal A}_b b^sr^s_0\frac{a}{b\gamma_4'}\inf\{\frac{k_b}{\#{\mathcal A}_b}:\ 0<b<1\}.
\end{array}
\]
Due to the above  fact (b), we see that $\sup\{\#{\mathcal A}_b b^s:\ 0<b<1\}<+\infty$. Hence, Theorem \ref{th4.1} (i) shows $\overline{\dim}_{\rm B}K_1\leq s$.
  Together with \eqref{eq4.15}, we  get $\dim_{H}K_1=s$, so $\dim_{H}K_1=\underline{\dim}_{\rm B}K_1=\overline{\dim}_{\rm B}K_1=s$. Finally, equality \eqref{eq2.2} shows $\dim_{H}K_n=\underline{\dim}_{\rm B}K_n=\overline{\dim}_{\rm B}K_n=s$ for all $n>0$. Statement (ii) is proven.
\end{pf}

\section{Remarks and Examples}
\setcounter{equation}{0}\setcounter{theorem}{0}

The following is an example of a generalized self-affine MIFS not satisfying the MBDP and the MWHP.

\begin{example}\label{Ex5.1}
Let $\Phi_n=\left\{\phi_0,\ \phi_1\right\}$ for all $n>0$, where
\[
\phi_0\left(
                  \begin{array}{c}
                    x\\
                    y \\
                  \end{array}
                \right)=\left(
                  \begin{array}{cc}
                    0.5 & 0 \\
                    0 & 0.4 \\
                  \end{array}
                \right)\left(
                  \begin{array}{c}
                    x\\
                    y \\
                  \end{array}
                \right),\quad
\phi_1\left(
                  \begin{array}{c}
                    x\\
                    y \\
                  \end{array}
                \right)=0.5\left(
                  \begin{array}{c}
                    x+1\\
                    y+1 \\
                  \end{array}
                \right).
\]
Then the MIFS does not satisfy the MWHP and the WBDP.
\end{example}
\begin{pf}
Let $I_n=0^n$ and $J_n=1^n$. It is clear that both and $I_n$ and $J_n$ belongs to $\mathcal{I}_{\frac{1}{2^n}}.$ A routine computation gives rise to
$$\frac{R_{I_n}}{r_{I_n}}=\left(\frac{5}{4}\right)^n \to +\infty,\ n\to+\infty,$$
which implies that the MBDP does not hold. Note that
\[
\phi_{J_n}(x)-\phi_{J_n}(y)=0.5^{n}(x-y),\quad \phi_{I_n}(x)-\phi_{{I_n}}(y)=\left(
                  \begin{array}{cc}
                    0.5^{n} & 0 \\
                    0 & 0.4^{n} \\
                  \end{array}
                \right)(x-y).
\]
By taking $x-y=(0,1)^t$, we have
\[
\frac{\|\phi_{J_{n}}(x)-\phi_{J_{n}}(y)\|}{\|\phi_{I_{n}}(x)-\phi_{I_{n}}(y)\|}=0.5^{n}\times 0.4^{-n}=\left(\frac54\right)^{n}.
\]
That means  the MIFS does not satisfy the MWHP.

\end{pf}

The following example shows that the conclusions of Corollary \ref{th4.4} (ii) and Theorem \ref{th4.5} do not hold if we omit the condition \eqref{eq3.2} in the definition of MOSC.
Specifically, without the condition \eqref{eq3.2} in the definition of MOSC, the conclusions $\underline{\dim}_{\rm B}K_n=\liminf\limits_{b\to0^+}\frac{\ln\#{\mathcal A}_b}{-\ln b}$ and $\limsup\limits_{b\to0^+}\frac{\ln\#{\mathcal A}_b}{-\ln b}=\overline{\dim}_{\rm B}K_n$ in Theorem \ref{th4.1} do not hold. Moreover, $\underline{\dim}_{\rm B}K_n$, $\overline{\dim}_{\rm B}K_n$ and $\dim_{H}(K_n)$ maybe independent of contraction ratios of the MIFS.

\begin{example}\label{Ex5.3}

Let the self-similar MIFS $\left\{\phi_{n,j}\right\}_{n\geq1}$ be defined on $\mathbb R$ with $0<\rho\le 1$ as
\[
\phi_{n,0}(x)=\frac12x, \quad \phi_{n,1}(x)=\frac12\left(x+\frac{1}{n}\rho^n\right).
\]
The following statements hold.
\medskip

{\rm (i)}  \eqref{eq3.1} holds with open sets $V_n=\left(0,\ \sum_{k=1}^\infty\frac{\rho^{n+k-1}}{2^{k}(n+k-1)}\right)$, and  \[K_{n}=\left\{\sum_{k=1}^\infty\frac{\rho^{n+k-1}j_k}{2^{k}(n+k-1)}:\ j_k=0,\ 1\right\}\subseteq\left[0,\ \sum_{k=1}^\infty\frac{\rho^{n+k-1}}{2^{k}(n+k-1)}\right]\] for  any $n\geq1$. But \eqref{eq3.2} does not hold for any open sets $V_n$ satisfying \eqref{eq3.1}.

\medskip

{\rm (ii)}  $\dim_{H}(K_n)=\underline{\dim}_{\rm B}(K_n)=\overline{\dim}_B(K_n)=\frac{\ln2}{\ln2-\ln\rho}$ and ${\mathcal H}^{\frac{\ln2}{\ln2-\ln\rho}}(K_n)=0$ for all $n>0$.

\medskip

{\rm (iii)}  Both \eqref{eq3.7} and the MWSC are not satisfied.
\end{example}
\begin{pf}
(i) It  is clear  that for any $n\geq1$, we have
\[
\sum_{k=1}^\infty\frac{\rho^{n+k-1}j_k}{2^{k}(n+k-1)}<\sum_{k=1}^\infty\frac{\rho^{n+k-1}}{2^{k}n}=\frac{\rho^n}{n(2-\rho)}.
\]
In combination with the above inequality, by the definition of the MIFS, a direct computation shows
\[\begin{array}{rl}
& K_{n}=\left\{\sum_{k=1}^\infty\frac{\rho^{n+k-1}j_k}{2^{k}(n+k-1)}:\ j_k=0,\ 1\right\}\medskip\\
\subseteq &[0,\ \sum_{k=1}^\infty\frac{\rho^{n+k-1}}{2^{k}(n+k-1)}]\subset [0,\ \frac{\rho^n}{n(2-\rho)}),\quad \forall\ n>0.
\end{array}\]
Hence, we have $V_n\subset (0,\ \frac{\rho^n}{n(2-\rho)})$,
\[
\phi_{n,0}(V_{n+1})\subset\left(0,\ \frac{\rho^{n+1}}{2(n+1)(2-\rho)}\right)\subset\left(0,\ \frac{\rho^n}{2n}\right)
\]
 and
\[
\phi_{n,1}(V_{n+1})=\left(\frac{\rho^n}{2n},\ \sum_{k=1}^\infty\frac{\rho^{n+k-1}}{2^{k}(n+k-1)}\right)
\]
for $n=1,2,\cdots$. That means  \eqref{eq3.1} holds.

\medskip

Based on observation that  $\phi_{n,0}^{-1}\circ\phi_{n,1}(x)=x+\frac{1}{n}\rho^n$ tends to the identity map, we conclude that \eqref{eq3.2} does not hold for any open sets $V_n$ satisfying \eqref{eq3.1}.

\medskip

{\rm (ii)}  We need only to consider the case $n=1$. Noting that
\[
K_{1}=\left\{\sum_{k=1}^\infty\frac{\rho^{k}j_k}{2^{k}k}:\ j_k=0,\ 1\right\},
\]
we see that $K_1$ is also generated by the MIFS $\left\{\psi_{n,j}(x)=r_{n}(x+j),\ j=0,\ 1\right\}_{n=1}^{+\infty}$ with $r_1=\frac{\rho}{2}$ and $r_n=\frac{(n-1)\rho}{2n}$ when $n\ge 2$. It is easy to check that the MIFS $\left\{\psi_{n,j},\ j=0,\ 1\right\}_{n=1}^{+\infty}$ satisfies the MOSC with the open set $V_n=(0,\frac{\rho}{2-\rho}).$   Hence, Theorem \ref{th4.3} shows $\underline{\dim}_B(K_1) =\overline{\dim}_B(K_1) =\dim_{H}(K_n)=\frac{\ln2}{\ln2-\ln\rho}$ and ${\mathcal H}^{\frac{\ln2}{\ln2-\ln\rho}}(K_n)=0$ by noting that the $s_k$ in Corollary  \ref{th4.4} equals $\frac{\ln2}{\ln2-\ln\rho+k^{-1}\ln k}$ for all $k>0$. This shows that the dimension of $K_1$ maybe any number in the interval $(0,1]$ when $\rho$ varies in the interval $(0,1]$.

\medskip

{\rm (iii)}  The proof of the conclusion (i) shows that $\phi_{n,0}^{-1}\circ\phi_{n,1}(x)=x+\frac{1}{n}\rho^n$ tends to the identity map. Let $X=[0,\ \sum_{k=1}^\infty\frac{\rho^{k}}{2^{k}k}]$. The self-similar MIFS $\left\{\phi_{n,j},\ j=0,\ 1\right\}_{n=1}^{+\infty}$ is well-defined on $X$. Let $I_n=0^{2n}$ and $J_{n,k}=0^k10^{2n-k-1}$.
By a simple calculation, we get 
\[
\phi_{1,I_n}^{-1}\circ\phi_{1,J_{n,k}}(x)=x+\frac{2^{2n-k-1}\rho^{k+1}}{k+1}<x+\frac{1}{k+1},\quad \frac{2n\ln2}{\ln2-\ln\rho}\le k< 2n.
\]
Then, we have
\[
\lim_{n\to\infty}\sup\left\{|\phi_{1,I_n}^{-1}\circ\phi_{1,J_{n,k}}(x)-x|:\ 2n> k\ge \frac{2n\ln2}{\ln2-\ln\rho} ,\ x\in X\right\}=0.
\]
Therefore,
\[
\#\left\{\phi_{1,J}:\ \phi_{1,I_n}(X)\cap\phi_{1,J}(X)\ne \emptyset\right\}\ge\#\left\{\phi_{1,J_{n,k}}:\ 2n> k\ge \frac{2n\ln2}{\ln2-\ln\rho}\right\}\to+\infty
\]
as $n\to \infty$. Hence, \eqref{eq3.7} does not hold.

\medskip

 From Lemma \ref{th3.5} (i), it follows that the MIFS does not satisfy the MWSC by noting that our MIFS satisfies the MWHP. In other words, there does not exist any sequence of subsets $U_n\subseteq X$ such that both \eqref{eq3.3} and \eqref{eq3.4} hold.
\end{pf}

Similar to Example \ref{Ex5.3}, we have the following example which shows that the Lebesgue measure in \eqref{eq3.2} and \eqref{eq3.3} can not be replaced by the diameter.

\begin{example}\label{Ex5.4}
Let the self-similar MIFS $\left\{\phi_{n,j}\right\}_{n\geq1}$ be defined on $\mathbb R^2$     as
$$\phi_{n,i+2j}((x,y)^t)=\left(\frac12\left(x+\frac{i}{2^nn}\right),\ \frac12\left(y+\frac{j}{2}\right)\right)^t,\ i,\ j=0,\ 1.$$
The following statements hold.

{\rm (i)}  \eqref{eq3.1} holds with open sets $$V_n= \left(0,\ \sum_{k=1}^\infty\frac{1}{2^{n+2k-1}(n+k-1)}\right)\otimes\left(0,\ \frac{1}{2}\right), \ \ n=1,2,\cdots$$
and
\[\begin{array}{rl}
K_{n}=&\left\{\sum_{k=1}^\infty(\frac{j_k}{2^{n+2k-1}(n+k-1)},\ \frac{i_k}{2^{k}})^t:\ i_k,\ j_k=0,\ 1\right\}\\
\subset& [0,\ \sum_{k=1}^\infty\frac{1}{^{n+2k-1}(n+k-1)}]\otimes[0,\ \frac{1}{2}],\ n=1,2,\cdots.
\end{array}
\]But \eqref{eq3.2} does not hold for any open sets $V_n$ satisfying \eqref{eq3.1}.

\medskip

{\rm (ii)}  $\dim_{H}(K_n)=\frac32$ and ${\mathcal H}^{\frac32}(K_n)=0$.

\medskip

{\rm (iii)}  Both \eqref{eq3.7} and MWSC are not satisfied.
\end{example}

The following example shows that, in Theorem \ref{th4.3}, it is possible that any one of ${\mathcal H}^s(K_n)=0$, ${\mathcal H}^s(K_n)=+\infty$ and $0<{\mathcal H}^s(K_n)<+\infty$  holds even if all assumptions in Theorem \ref{th4.3} have been satisfied.

\begin{example}\label{Ex5.5}
Let $a_n\in [\frac23,\ \frac32]$ be real numbers such that $\lim\limits_{n\to\infty}3^{nt}\prod\limits_{j=1}^na_j=0$ for all $t<0$ and $\lim\limits_{n\to\infty}3^{nt}\prod\limits_{j=1}^na_j=+\infty$ for all $t>0$.
Define the self-similar MIFS $\left\{\phi_{n,j}\right\}_{n\geq1}$ on $\mathbb R$ as
\[\phi_{n,0}(x)=\frac1{3a_n}x,\ \ \phi_{n,1}(x)=\frac1{3a_n}(x-1)+1.\]
The following statements hold.
\medskip

{\rm (i)}  The MIFS satisfies the MOSC defined in Definition \ref{def3.1} with open sets $V_n=(0,\ 1)$ for all $n\geq1$. Furthermore, we have $\{0,\ 1\}\subset K_{n}\subset[0,\ 1]$ for all $n\geq1$.

\medskip

{\rm (ii)}  $\dim_{H}(K_n)=\frac{\ln2}{\ln3}$.

\medskip

{\rm (iii)}  ${\mathcal H}^{\frac{\ln2}{\ln3}}(K_n)=0$ when $\limsup\limits_{n\to\infty}\prod\limits_{j=1}^na_j=+\infty$.

\medskip

{\rm (iv)}  ${\mathcal H}^{\frac{\ln2}{\ln3}}(K_n)=+\infty$ when $\lim\limits_{n\to\infty}\prod\limits_{j=1}^na_j=0$.

\medskip

{\rm (v)}  $0<{\mathcal H}^{\frac{\ln2}{\ln3}}(K_n)<+\infty$ when $0<\limsup\limits_{n\to\infty}\prod\limits_{j=1}^na_j<+\infty$.
\end{example}

{\bf Remark}. (a) It is easy to see  that for any integer $j\geq1$, we have  $\frac{j+3}{j+2}\in (1,\ \frac32)$, $\limsup\limits_{n\to\infty}\prod\limits_{j=1}^n(\frac{j+3}{j+2})$ $=+\infty$ and $\lim\limits_{n\to\infty}3^{nt}\prod\limits_{j=1}^n(\frac{j+3}{j+2})=0$ for all $t<0$.
(b) It is also obvious that $\frac{j+2}{j+3}\in (\frac23,\ 1)$ satisfy $\lim\limits_{n\to\infty}\prod\limits_{j=1}^n(\frac{j+2}{j+3})=0$ and $\lim\limits_{n\to\infty}3^{nt}\prod\limits_{j=1}^n(\frac{j+2}{j+3})=+\infty$ for all $t>0$. (c) One can check that $(\frac{3}{2})^{2^{-j}}\in (1,\ \frac32)$ and $\limsup\limits_{n\to\infty}\prod\limits_{j=1}^n(\frac{3}{2})^{2^{-j}}=\frac{3}{2}$. Thus,
$\lim\limits_{n\to\infty}3^{nt}\prod\limits_{j=1}^n(\frac{3}{2})^{2^{-j}}$ equals $+\infty$ or $0$ based on $t$ being positive or negative, respectively.

\begin{pf}
(i) Actually, by the definitions of the MIFS, it is easy to see that $\phi_{n,0}(0)=0$, $\phi_{n,1}(1)=1$ and $\phi_{n,j}([0,\ 1])\subset [0,\ 1]$ for all $n\ge1$ and $j=0,\ 1$. Hence $\{0,\ 1\}\subset K_{n}\subset[0,\ 1]$ for all $n\ge1$.

\medskip

By the assumption $a_n\in [\frac23,\ \frac32]$, we see that $\phi_{n,0}((0,\ 1))\cap\phi_{n,1}((0,\ 1))=\emptyset$ for all $n\ge1$, so the MOSC holds. The statement (i) follows.

\medskip

{\rm (ii)} By the definition of our MIFS, we have
\[
\prod\limits_{i=1}^n\sum\limits_{j=1}^{2}r_{i,j}^s=2^n\left(3^{n}\prod\limits_{j=1}^na_j\right)^{-s},\quad \forall \ s\in{\mathbb R}.
\]
Since $\lim\limits_{n\to\infty}3^{nt}\prod\limits_{j=1}^na_j$ is 0 when $t<0$ and $+\infty$ when $t>0$, we see that  $\dim_{H}(K_1)=\frac{\ln2}{\ln3}$ by using Theorem \ref{th4.3} or Corollary \ref{th4.4}.  Using \eqref{eq2.2}, one see that $\dim_{H}(K_n)=\dim_{H}(K_1)$. We finish the proof of  the statement (ii).

\medskip

{\rm (iii)}  By letting $s=\frac{\ln2}{\ln3}$ and $b_n=\left(3^{n}\prod\limits_{j=1}^na_j\right)^{-1}$, we have
\begin{equation}\label{eq5.3}
{\mathcal H}^s_{b_n}(K_1)\le \sum_{I\in {\mathcal I}_{b_n}}|\phi_{1,I }([0,\ 1])|^s=2^n\left(3^{n}\prod\limits_{j=1}^na_j\right)^{-s}=\left(\prod\limits_{j=1}^na_j\right)^{-s}.
\end{equation}
Since ${\mathcal H}^s(K_1)=\lim\limits_{\delta\downarrow 0}{\mathcal H}^s_{\delta}(K_1)$ and ${\mathcal H}^s_{\delta}(K_1)$ is a decreasing function of $\delta$ in interval $(0,\ 1)$, we see that ${\mathcal H}^s(K_1)=\lim\limits_{n\to\infty}{\mathcal H}^s_{b_n}(K_1)\le\liminf\limits_{n\to\infty}\left(\prod\limits_{j=1}^na_j\right)^{-s}=0$ by using the assumption $\limsup\limits_{n\to\infty}\prod\limits_{j=1}^na_j=+\infty$. Therefore, ${\mathcal H}^s(K_1)=0$. According to  \eqref{eq2.2}, one see that $\dim_{H}(K_n)=0$ for all $n>0$. Hence, the statement (iii) is proven.

\medskip

{\rm (iv)}
As we have done in the proof of (iii), write $b_n=\left(3^{n}\prod\limits_{j=1}^na_j\right)^{-1}$ for all $n>0$. Let $\{\mu_n\}_{n\geq1}$ be defined by \eqref{eq2.3} with our MIFS and the probability weights $p_{n,j}=\frac12$ ($n\ge1,\ j=1,\ 2$).

\medskip

For any small number $\delta>0$ and set $U\subset{\mathbb R}$ satisfying $|U|\le \delta$, there is an integer $n>0$ such that $b_{n+1}\le |U|<b_{n}$ by noting that $b_n$ decreases to 0 as $n$ tends to infinity. It follows that $U\subseteq(a,\ a+b_n)$ for some $a\in{\mathbb R}$.
For any two distinct words $I\ne J\in\Sigma_1^n$, statement (i) shows that $\phi_{1,I}([0,\ 1])$ and $\phi_{1,J}([0,\ 1])$ are contained in two disjoint intervals with length $b_n$. Therefore, $U$ intersects at most two subsets of $\{\phi_{1,I}([0,\ 1]):\ I\in\Sigma_1^n\}\}$.
Hence, we have
\begin{equation}\label{eq5.4}
\mu_1(U)\le\frac2{2^n}=\frac2{3^{ns}}=2\left(3\prod\limits_{j=1}^{n+1}a_j\right)^sb_{n+1}^s \le 2\left(3\prod\limits_{j=1}^{n+1}a_j\right)^s|U|^s,
\end{equation}
where $s=\frac{\ln2}{\ln3}$ is the Hausdorff dimension of $K_n$. Using the definition of Hausdorff measures shows
\[\begin{array}{rl}
{\mathcal H}^s(K_1)=&\lim\limits_{\delta\to0}\inf\{\sum\limits_{j}|U_j|^s:\ |U_j|\le\delta,\ \bigcup\limits_{j}U_j\supseteq K_1\}\medskip\\
\geq&\lim\limits_{n\to\infty}\frac12\left(3\prod\limits_{j=1}^{n+1}a_j\right)^{-s}=+\infty
\end{array}\]
by using $\lim\limits_{n\to\infty}\prod\limits_{j=1}^na_j=0$. Using \eqref{eq2.2}, one see that $\dim_{H}(K_n)=+\infty$ for any $n\geq2$. Thus, we prove  the statement (iv).

\medskip

{\rm (v)}  The inequality \eqref{eq5.3} shows that ${\mathcal H}^{\frac{\ln2}{\ln3}}(K_n)<+\infty$ when $0<\limsup\limits_{n\to\infty}\prod\limits_{j=1}^na_j$. On the other hand, by the assumption $\limsup\limits_{n\to\infty}\prod\limits_{j=1}^na_j<+\infty$, we get  $$c:=\sup\left\{2\left(3\prod\limits_{j=1}^{n+1}a_j\right)^{\frac{\ln2}{\ln3}}:\ n>0\right\}<+\infty.$$ From \eqref{eq5.4}, it follows  that $\mu_1(U)\le c |U|^{\frac{\ln2}{\ln3}}$ for all subset $U\subset{\mathbb R}$. By the well known mass distribution principle, we get that ${\mathcal H}^{\frac{\ln2}{\ln3}}(K_1)\ge \frac1c>0$. It is   obviously true that ${\mathcal H}^{\frac{\ln2}{\ln3}}(K_n)>0$ by using \eqref{eq2.2}.  The proof of statement (v) is completed.
\end{pf}

The following example shows that \eqref{eq3.7} maybe wrong if we remove the MWHP in Lemma \ref{th3.5}.

\begin{example}\label{Ex5.6}
Consider the self-similar MIFS $\left\{\ \phi_{n,1},\ \cdots,\ \phi_{n,n^2+3}\right\}_{n=1}^{+\infty}$ on ${\mathbb R}^2$ with
\[\left\{\begin{array}{ll}
\phi_{n,i+1+nj}((x,\ y)^t)=\frac1{2n}((x,\ y)^t+(i,\ j)^t),\quad& 0\le i,\ j\le n-1,\\
\phi_{n,i+n^2+2j}((x,\ y)^t)=\frac1{2}((x,\ y)^t+(i,\ j)^t),\quad& \ 0\le i,\ j\le 1,\ i+j>0.
\end{array}\right.
\]
Then, the MIFS satisfies the MOSC defined in Definition \ref{def3.1}, but it does not satisfy the MWHP defined in Definition \ref{def3.3}. Furthermore, the inequality \eqref{eq3.4} does not hold when we take $U_n=[0,\ 1]^2$, and the inequality \eqref{eq3.7} does not hold when we take $X=[0,\ 1]^2$.
\end{example}
\begin{pf}
Let $U=(0,\ 1)^2$. It is easy to check that $\phi_{n,i+1+jn}(U)=(\frac{i}{2n},\ \frac{i+1}{2n})\times(\frac{j}{2n},\ \frac{j+1}{2n})$ when $0\le i,\ j\le n-1$ and, $\phi_{n,i+n^2+2j}(U)=(\frac{i}{2},\ \frac{i+1}{2})\times(\frac{j}{2},\ \frac{j+1}{2})$ when $0\le i,\ j\le 1$ and $i+j>0$. They are $n^2+3$ disjoint subsets of $U$. That means that  the MIFS satisfies the MOSC defined in Definition \ref{def3.1} with open sets $V_n=U$ .

\medskip

It is easy to see that  $K_n=[0,\ 1]^2$ for all $n\ge1$. Let $X=[0,\ 1]^2$. We see that  $\left\{\ \phi_{n,1},\ \cdots,\ \phi_{n,n^2+3}\right\}_{n=1}^{+\infty}$ is an MIFS on $X$. Take $b_n=\frac{1}{2^{n}(n-1)!}$ and  $I_{n-1}=\underbrace{11\cdots 1}_{n-1}$. From a direct computation, it follows  that the contraction ratio of $\phi_{1,I_{n-1}j}$ is $\frac{1}{2^nn!}$ if $1\le j\le n^2$ and $\frac{1}{2^n(n-1)!}$ if $n^2<j\le n^2+3$. This implies that $\phi_{1,I_{n-1}j}$ belongs to ${\mathcal A}_{b_n}$ for all $1\le j\le n^2+3$. Therefore, the MIFS  does not satisfy the MWHP since the supremum in \eqref{eq3.5} is $+\infty$.

\medskip

It is easy to check that
\[
 \phi_{n,n^2+1}(X)\cap\phi_{n,n+nj_n}(X)\neq\emptyset,\quad 0\le j_n\le n-1.
\]
Then, we have
\[
  \phi_{1,I_{n-1}}\circ\phi_{n,n^2+1}(X)\cap\phi_{1,I_{n-1}}\circ\phi_{n,n+nj_n}(X)\neq\emptyset,\quad 0\le j_n\le n-1.
\]
Hence, we conclude that the \eqref{eq3.4} does not hold when we take $U_n=[0,\ 1]^2$, since the supremum in \eqref{eq3.4} is $+\infty$ in this case. The same argument shows that the inequality \eqref{eq3.7} does not hold when we take $X=[0,\ 1]^2$.
\end{pf}

The following example shows that it is possible that $\underline{\dim}_{\rm B}K_n<\overline{\dim}_{\rm B}K_n$  even if all assumptions in Theorem \ref{th4.3} are satisfied.

\begin{example}\label{Ex5.7}
Consider the self-similar MIFS $\left\{\phi_{n,0},\ \phi_{n,1}\right\}_{n=1}^{+\infty}$ on $\mathbb R$ with
\begin{equation}\label{eq5.6}\left\{\begin{array}{ll}
\phi_{n,j}(x)=\frac1{3}(x-j)+j,\quad& 2^{2k}\le n<2^{2k+1},\ j=0,\ 1,\ k\ge0,\\
\phi_{n,j}(x)=\frac1{2}(x-j)+j,\quad& 2^{2k+1}\le n<2^{2k+2},\ j=0,\ 1,\ k\ge0.
\end{array}\right.\end{equation}
Then, the MIFS satisfies the MOSC, MWSC, MBDP and MWHP. Furthermore, we have $\dim_{H}K_n=\underline{\dim}_{\rm B}K_n=\frac{3\ln2}{2 ln3+\ln2}<\overline{\dim}_{\rm B}K_n=\frac{3\ln2}{\ln3+2\ln2}$.
\end{example}
\begin{pf}
Actually, by the definition of the MIFS, it is easy  to see $\phi_{n,0}(0)=0$, $\phi_{n,1}(1)=1$, $\phi_{n,j}([0,\ 1])\subset [0,\ 1]$ and $\phi_{n,0}((0,\ 1))\cap \phi_{n,1}((0,\ 1))=\emptyset$
for $n\geq1$ and $j=0,\ 1$. Hence, the MIFS satisfies the MOSC, MWSC, MBDP and MWHP.

\medskip

That means for any $n\geq1$, we have $\mathcal{A}_{r_1\cdots r_n}=\{\phi_{1,I}:I\in\Sigma_1^n\}$, which implies that $\# \mathcal{A}_{r_1\cdots r_n}=2^n.$
From \eqref{eq5.6}, a direct computation leads to
\begin{equation*}
\frac{\ln \#\mathcal{A}_{r_1\cdots r_n}}{-\ln(r_1\cdots r_n)}=
\left\{
\begin{array}{l}
\frac{n\ln2}{\left(\frac{2}{3}(1-4^k)+n\right)\ln3+\frac{2}{3}(4^k-1)},\ \quad 2^{2k}\leq n<2^{2k+1},\\
\frac{n\ln2}{\frac{1}{3}(4^{k+1}-1)\ln3 +(\frac{1}{3}(1-4^{k+1})+n)\ln2}, \ \ 2^{2k+1}\leq n<2^{2k+2}.
\end{array}\right.
\end{equation*}
By noting that \eqref{contra-ra} holds, from  \eqref{lower-dim}, it follows that
\[\underline{\dim}_{\rm B}K_1=\liminf\limits_{n\to\infty}\frac{\ln \#\mathcal{A}_{r_1\cdots r_n}}{-\ln(r_1\cdots r_n)}=\frac{3\ln2}{2 \ln3+\ln2}.\]
Also, by \eqref{upper-dim},  we have
\[\overline{\dim}_{\rm B}K_1=\limsup\limits_{n\to\infty}\frac{\ln \#\mathcal{A}_{r_1\cdots r_n}}{-\ln(r_1\cdots r_n)}=\frac{3\ln2}{ \ln3+2\ln2}.\]

\medskip

Finally, it is easy to see that the $s_k$ defined in Corollary \ref{th4.4} equals $\frac{\ln2^{k}}{-\ln R_I}=\frac{\ln2^{k}}{-\ln r_I}$ with $I\in\Sigma_1^k$. Hence $\liminf\limits_{k\to\infty}s_k=\underline{\dim}_{\rm B}K_1$ by using Corollary \ref{th4.4} (i). The proof is completed.
\end{pf}

The following example shows that the conclusion $\overline{\dim}_{\rm B}K_n<\limsup\limits_{b\to0^+}\frac{\ln\#{\mathcal A}_b}{-\ln b}$ is possible if we remove the assumption $r_0=\inf\{r_{n,j}:\ n\ge 1,\ 1\le j\le N_n\}>0$ in Theorem \ref{th4.1} even if we assume the MOSC and $r_{n,j}=r_n$. Moreover, the conclusion $\overline{\dim}_{\rm B}K_n<\limsup\limits_{b\to0^+}\frac{\ln\#{\mathcal A}_b}{-\ln b}$ depends on both the assumption $r_0>0$ and the translations in the maps of the MIFS.

\begin{example}\label{Ex5.8}Let $\left\{\Phi_n\right\}_{n=1}^\infty$ be a self-similar MIFS with $\Phi_n=\{\phi_{n,j}(x)=2^{-2^n}(x+\alpha_j):\ \alpha_j\in D_n\}$ and $\{0,\ 2^{2^n}-1\}\subseteq D_n\subseteq \{j:\ 0\leq j\leq 2^{2^n}-1\}$.

\medskip

{\rm (i)}  The MIFS satisfies the MOSC defined in Definition \ref{def3.1} with $V_n=(0,\ 1)$. Moreover, the left side of \eqref{eq3.7} is less than or equal to 3 if we choose $X=[0,\ 1]$.

\medskip

{\rm (ii)}  If $D_n=\{j:\ 0\leq j\leq 2^{2^n}-1\}$ for all $n\ge1$, then
$1=\dim_{\rm H}K_n=\dim_{\rm B}K_n=\dim_{\rm P}K_n=\liminf\limits_{b\to0^+}\frac{\ln\#{\mathcal A}_b}{-\ln b}<\limsup\limits_{b\to0^+}\frac{\ln\#{\mathcal A}_b}{-\ln b}$
for all $n\ge1$.

\medskip

{\rm (iii)}  If $D_n=\{0,\ 2^{2^n}-1\}$ for all $n\ge1$, then $0=\dim_{\rm H}K_n=\dim_{\rm B}K_n=\dim_{\rm P}K_n=\liminf\limits_{b\to0^+}\frac{\ln\#{\mathcal A}_b}{-\ln b}=\limsup\limits_{b\to0^+}\frac{\ln\#{\mathcal A}_b}{-\ln b}$
for all $n\ge1$.
\end{example}
\begin{pf}
(i)  It is clear that  $K_n=[0,\ 1]$ for all $n>0$, so $\phi_{n,j}([0,\ 1])\subset [0,\ 1]$ for all $n,\ j\geq1$, we can choose $X=[0,\ 1]$. By checking  $\phi_{n,i}((0,\ 1)))\cap\phi_{n,j}((0,\ 1))=\emptyset$ for all $n>0$ and $i\neq j$, we see that the MIFS satisfies the MOSC defined in Definition \ref{def3.1}. The second conclusion is obvious.

\medskip

{\rm (ii)}  Since $K_n=[0,\ 1]$ for all $n>0$, the definition of our MIFS shows
\begin{equation}\label{eq5.7}
R_{1,j_1j_2\cdots j_n}=r_{1,j_1j_2\cdots j_n}=\prod\limits_{i=1}^n2^{-2^i}=2^{(-2^{n+1}+2)},\ \ n=1,2,\cdots.
\end{equation}

Let $N_\delta$ be the number of $\delta$-mesh cubes that intersect $K_1=[0,\ 1]$.
By the definition of our MIFS, every interval $\phi_{1,j_1j_2\cdots j_{n}}([0,\ 1])$ is exactly a $\delta$-mesh cube when $\delta=2^{(-2^{n+1}+2)}$. Hence, we have
\begin{equation}\label{eq5.8}
N_{2^{(-2^{n+1}+2)}}=\#{\mathcal A}_{2^{(-2^{n+1}+2)}}=2^{(2^{n+1}-2)}.
\end{equation}
The definition of ${\mathcal A}_b$ and \eqref{eq5.7} show
$\#{\mathcal A}_b=2^{(2^{n+1}-2)}$ when $2^{(-2^{n+1}+2)}\le b< 2^{(-2^{n}+2)}$. Consequently, it follows that
\begin{equation}\label{eq5.9}
\sup_{2^{(-2^{n+1}+2)}\le b< 2^{(-2^{n}+2)}}\frac{\ln(\#{\mathcal A}_b)}{-\ln b}=\frac{2^{n+1}-2}{2^{n}-2}
\end{equation}
and
\[
\inf_{2^{(-2^{n+1}+2)}\le b< 2^{(-2^{n}+2)}}\frac{\ln(\#{\mathcal A}_b)}{-\ln b}=1.
\]
Therefore
\[
\liminf\limits_{b\to0^+}\frac{\ln\#{\mathcal A}_b}{-\ln b}=1,\ \ \limsup\limits_{b\to0^+}\frac{\ln\#{\mathcal A}_b}{-\ln b}=2.
\]
Hence statement (ii) follows by using $K_n=[0,\ 1]$.

\medskip

{\rm (iii)}  By a similar argument to \eqref{eq5.8}, we have $N_{2^{(-2^{n+1}+2)}}=\#{\mathcal A}_{2^{(-2^{n+1}+2)}}=2^n$. Hence, we get  $N_b\leq 2N_{2^{(-2^{n+1}+2)}}=2^{n+1}$ when $2^{(-2^{n+1}+2)}\le b< 2^{(-2^{n}+2)}$. Also, by a similar argument to \eqref{eq5.9}, we have
\[
\sup_{2^{(-2^{n+1}+2)}\le b< 2^{(-2^{n}+2)}}\frac{\ln(\#{\mathcal A}_b)}{-\ln b}=\frac{n}{2^{n}-2}.
\]
Therefore, $\limsup\limits_{b\to0^+}\frac{\ln N_b}{-\ln b}=\limsup\limits_{b\to0^+}\frac{\ln\#{\mathcal A}_b}{-\ln b}=0$, and so our conclusion (iii) follows.
\end{pf}

The following example shows that the $U_n$ in \eqref{eq3.4} can not be chosen arbitrarily in general cases. The inequalities   \eqref{eq3.7} and \eqref{eq3.8} depend on the basic space $X$. On the other hand, the MWHP is not necessary condition for \eqref{eq3.7} and \eqref{eq3.8}.

\begin{example}\label{Ex5.9}
Consider the self-similar MIFS $\left\{\phi_{n,0},\ \phi_{n,1},\ \cdots,\ \phi_{n,n}\right\}_{n=1}^{+\infty}$ on $X\subset\mathbb R$ with
\[
\phi_{n,0}(x)=\frac1{2}(x+1),\ \
\phi_{n,j}(x)=\frac1{2n}(x+j-1),\quad 1\le j\le n,\ n=1,\ 2,\ \dotsc.
\]
Then, the MIFS satisfies the MOSC but does not satisfy the MWHP.  
The inequality \eqref{eq3.4} does not hold when $U_n^o\supset[0,\ 1]$, \eqref{eq3.7} and \eqref{eq3.8} do not hold when $X^o\supset[0,\ 1]$.
\end{example}
\begin{pf}
A direct computation shows $K_n=[0,1]$ and $\phi_{n,i}((0,1))\cap \phi_{n,j}((0,1))=\emptyset$ for all $i\ne j\in\{0,\ 1,\ \cdots,\ n\}$. It follows that the MIFS satisfies the MOSC. By a similar argument to Example \ref{Ex5.6}, we see the MIFS has no MWHP.

Let $\epsilon>0$  such that $U\supset[-\epsilon,\ 1+\epsilon]$. Choose $b_n=\frac{19}{2^n10}$ ($n>2$). It is clear that $0^{n-1}j\in{\mathcal I}_{b_n}$ for any $0\le j\le n$. By a simple argument, we conclude that for $n\ge j\ge n(1-\epsilon)-\epsilon$,
\[
\phi_{n,0}(U)\cap\phi_{n,j}(U)\supset \left[\frac{1-\epsilon}{2},\ \frac{2+\epsilon}{2}\right]\bigcap
\left[\frac{j-1-\epsilon}{2n},\ \frac{\epsilon+j}{2n}\right]\ne\emptyset.
\]
By noting that \[
\phi_{1,0^{n}}(U)\cap\phi_{1,0^{n-1}j}(U)=\phi_{1,0^{n-1}}(\phi_{n,0}(U)\cap\phi_{n,j}(U)),
\]
we  obtain that
\[
\#\{J\in{\mathcal I}_{b_n}:\ \phi_{1,0^{n}}(U)\cap\phi_{1,0^{n-1}j}(U)\ne\emptyset\}\ge n\epsilon\to+\infty\mbox{ as }n\to+\infty.
\]
Therefore, \eqref{eq3.4} does not hold when we choose $U_n$ such that $U_n^o\supset[0,\ 1]$.

\medskip

The above proof also shows that \eqref{eq3.7} and \eqref{eq3.8} do not hold when $X^o\supset[0,\ 1]$.
However, it is easy to check that \eqref{eq3.7} and \eqref{eq3.8} hold when $X=[0,1]$.

\end{pf}

\bigskip

\small

\noindent School of Mathematics and Statistics, Key Laboratory of Analytical Mathematics and Applications (Ministry of Education), Fujian Key Laboratory of Analytical Mathematics and Applications (FJKLAMA), Fujian Normal University, 350117 Fuzhou, P.R. China.

\bigskip

E-mail: 1784336267@qq.com (Y.-S. Cao), dengfractal@126.com (Q.-R. Deng),

limtwd@fjnu.edu.cn~(M.-T. Li)

\end{document}